\documentclass[]{article}

\usepackage[letterpaper, margin=1in]{geometry}

\usepackage{ifpdf}

\usepackage{cite}

\usepackage{amssymb,latexsym,float,epsfig,algorithm,subcaption,algorithmicx,algpseudocode,tabularx,pdflscape,epigraph,xparse,color,cite,url}
\usepackage[cmex10]{amsmath}
\usepackage{fixltx2e}
\usepackage{url}
\newtheorem{theorem}{Theorem}[section]
\newtheorem{lemma}[theorem]{Lemma}
\begin{document}
\title{Surrogate Lagrangians for Variational Integrators: \\ High Order Convergence with Low Order Schemes}

\author{Gerardo~De~La~Torre~and~Todd.~D.~Murphey
\thanks{G. De La Torre and T. D. Murphey are with the Department of Mechanical Engineering, Northwestern University, Evanston, IL, 60208 USA  e-mail: \textit{gerardo.delatorre@northwestern.edu, t-murphey@northwestern.edu}.}
}

\maketitle

\begin{abstract}
Variational integrators are momentum-preserving and symplectic numerical methods used to propagate the evolution of Hamiltonian systems. 
In this paper, we introduce a new class of variational integrators that achieve fourth-order convergence despite having the same integration scheme as traditional second-order variational integrators.
The new class of integrators are created by replacing a dynamical system's Lagrangian in the variational integration algorithm with its surrogate Lagrangian. 
By incorporating the surrogate Lagrangian the propagation errors induced by variational integrators, up to a given order, are eliminated. 
Furthermore, no assumption on the Lagrangian's structure is made and, therefore, the proposed approach is applicable to a large range of dynamical systems.
In addition, surrogate variational integrators are also constructed for Hamiltonian systems subjected to holonomic constraints and external forces.  
Finally, the methodology is extended to derive higher-order surrogate variational integrators that achieve an arbitrary order of accuracy but retain second-order complexity in the integration scheme. 
Several numerical experiments are presented to demonstrate the efficacy of our approach.
\end{abstract}

\section{Introduction}
Variational integrators are well-suited for the time propagation of Euler-Lagrange equations since they arise through direct discretization of Hamilton's variational principle. 
As a result, variational integrators are able to ensure (or strongly enforce) the conservation of fundamental mechanical quantities such as momentum and energy \cite{West:thesis,VI:Marsdem:01, VI:book:01,VI:001}.
Furthermore, holonomic constraints, external forces, impacts, and non-smooth phenomenon fit naturally into the variational integration scheme \cite{fetecau2003nonsmooth}. 
From an implementation standpoint, variational integrators are scalable and implementable for generic mechanical systems in generalized coordinates \cite{VI:Murphey:01}. 
Furthermore, iterative projection-based optimization methods become real-time implementable when variational integrators are utilized despite relatively low sensor or actuator bandwidth  \cite{projection003, DLT:DISS}. 

The accuracy of traditional variational integrators is governed by the approximation used to define the discrete Lagrangian.    
Recently proposed variational integrators achieve an increase level of accuracy by using a variety of methods (e.g. Hermite interpolation, Galerkin methods, etc.) to obtain higher-order approximations of the discrete Lagrangian \cite{ober2015construction,leok2012general,Leok15112011,VI:Marsdem:01}.
However, these methods do not exploit the geometric properties of Hamiltonian systems nor those of variational integrators. 
Furthermore, the variational integration scheme is made more complex by introducing such approximations. 

Recently, backward error analysis has been used to proposed new integration methods.
Backward error analysis is used to quantify the modification to the propagated system induced by a particular integration scheme \cite{reich1999backward,hairer1994backward,WARMING1974159,griffiths1986scope}.
Numerical integrators are then constructed to mitigate the effects of the known modification \cite{vilmart2008study,chartier2007numerical}. 
When backward error analysis is applied to variational integrators the modification to the propagated system is described by a modified Hamiltonian or Lagrangian \cite{vermeeren2015modified,sanz1992symplectic}.
Therefore, variational integrators exactly capture the evolution of a ``near-by" Hamiltonian system. 
This property can be exploited by altering the integration scheme such that the ``near-by" system better represents the considered system.   
Using this idea numerical schemes that propagate the rotation of a rigid body and the evolution of mechanical systems with separable Hamiltonians have been proposed  \cite{chartier2007numerical,mushtaq2014higher,Mushtaq20141461}.

The main contribution of this paper is the presentation of a new class of variational integrators that achieve fourth-order convergence despite having the same integration scheme as traditional second-order variational integrators.
The new class of integrators are created by replacing a dynamical system's Lagrangian in the variational integration algorithm with its surrogate Lagrangian. 
By incorporating the surrogate Lagrangian the propagation errors induced by variational integrators, up to a given order, are eliminated. 
Furthermore, no assumption on the Lagrangian's structure is made and, therefore, the proposed approach is applicable to a large range of dynamical systems.
In addition, surrogate variational integrators are also constructed for Hamiltonian systems subjected to holonomic constraints and external forces.  
Finally, the methodology is extended to derive higher-order surrogate variational integrators that achieve an arbitrary order of accuracy but retain second-order complexity in the integration scheme. 
Several numerical experiments are presented to demonstrate the efficacy of our approach. 

The organization of this paper is as follows. 
Section 2 gives an overview of the midpoint variational integrator.  
Backward error analysis is reviewed in Section 3.
In addition, Section 3 formulates the modified Lagrangian for a conservative Hamiltonian system. 
Surrogate Lagrangians for conservative, forced, and constrained Hamiltonian systems are introduced in Section 4. 
Section 5 extends the presented methodology to produce higher order surrogate Lagrangians. 
Section 6 presents results from numerical experiments. 
Conclusions are discussed in Section~7.

\section{Variational Integrators}\label{sec:VI}
To begin our discussion, the formulation of the Euler-Lagrangian equations is reviewed by considering the sufficiently differentiable Lagrangian of a dynamical system represented as
\begin{align}
L(q(t), \dot q(t)) = T(q(t), \dot q(t)) - V(q(t)),
\end{align}
where $q$ is the state configuration vector, $\dot q$ is its time derivative, $T(q(t), \dot q(t))$ describes the system's kinetic energy, and $V(q(t))$ describes the system's potential energy.
The action, $S$, is defined as
 \begin{align} \label{back:action}
S[q(t)] = \int_{t_0}^{t_\textrm{f}}L(q(\tau), \dot q(\tau)) ~ \textrm{d} \tau.
\end{align}   
The least action principle is used to derive the variational relation 
\begin{align}
\delta S[q(t)] = \delta \int_{t_0}^{t_\textrm{f}}L(q(\tau), \dot q(\tau)) ~ \textrm{d} \tau = 0, \quad \delta q(t_0) = 0, \quad \delta q(t_\textrm{f})=0,
\end{align}
which results  in the classical Euler-Lagrange equations \cite{2002analytical}:
 \begin{align}\label{back:contEL}
\frac{\partial}{\partial  t}\frac{\partial L}{\partial  \dot q}(q, \dot q) - \frac{\partial L}{\partial  q}(q, \dot q) = 0.
\end{align}   
Propagating equation (\ref{back:contEL}) with numerical integration schemes developed for general second order differential equations will result in numerical errors since the system's fundamental characteristics (e.g. symmetries of motion, conservation of energy) are ignored. 
Variational integrators approximate the continuous trajectory of mechanical systems with a sequence of discrete points while ensuring (or strongly enforcing) the conservation of fundamental quantities such as momentum and energy \cite{VI:Marsdem:01}.  
Specifically, a sequence of system configuration vectors $\{(t_0,q_0), (t_1,q_1), \dots, (t_n,q_n)\}$ is found such that the continuous system trajectory is approximated as $q_m \approx q(t_m)$ where $h = t_{i+1}-t_i$ is the discretization time step. 
Derivations for the same variational integrator presented here are given in \cite{VI:Marsdem:01, ober2011discrete}.

The derivation of the \emph{midpoint} variational integrator begins by defining the \emph{discrete Lagrangian}, $L_\textrm{d}(q_k,q_{k+1})$, as
 \begin{align}\label{back:DL}
L_\textrm{d}(q_k,q_{k+1}) &= L(\frac{q_k + q_k}{2},\frac{q_{k+1} - q_k}{h})h \approx \int_{t_k}^{t_{k+1}}L(q(t),\dot{q}(t)) ~ h.
\end{align}   
A generalized midpoint approximation can be used to define other discrete Lagrangians. 
However, in this paper we consider the midpoint approximation since it results in second order accuracy as discussed in \cite{West:thesis} and shown in later in Section \ref{sec:backward:VI}.
Equation (\ref{back:action}) can be approximated as a sum of discrete Lagrangians:
\begin{align}
S[q(t)] \approx \sum_{k=0}^{n-1} L_\textrm{d}(q_k,q_{k+1}).
\end{align}   
It follows from the least action principle that
\begin{align}\label{back:actionsum}
\delta S \approx \sum_{k=1}^{n-1} (D_1L_\textrm{d}(q_k,q_{k+1}) + D_2L_\textrm{d}(q_{k-1},q_{k}))\cdot \delta q_{k} = 0,
\end{align}   
assuming $\delta q_0 = \delta q_n = 0$. 
The variations of the action sum are zero for any $\delta q_k$ and, as a result, the Discrete Euler-Lagrange (DEL) equation is derived as
\begin{align}
D_1L_\textrm{d}(q_k,q_{k+1}) + D_2L_\textrm{d}(q_{k-1},q_{k}) = 0.
\end{align}   
Notice that the DEL equation is the discrete time equivalent to the classical Euler-Lagrange equation (\ref{back:contEL}). 
Equivalently, the resulting DEL equation can be given its \emph{position-momentum} form as   
\begin{align}
&p_{k} + D_1L_\textrm{d}(q_k,q_{k+1}) = 0, \label{back:pm1}\\
&p_{k} =  D_2L_\textrm{d}(q_{k-1},q_{k}). \label{back:pm2}
\end{align}
Note that $p_{k}$ does not depend on $q_{k+1}$ and (in the unforced case) $p_{k}$ is the momentum quantity conserved by the integrator \cite{VI:Murphey:01, West:thesis}.
Furthermore, the previously defined two-step mapping $(q_{k-1},q_k)\rightarrow(q_{k+1})$ is now replaced with a one step mapping $(q_{k},p_k)\rightarrow(q_{k+1},p_{k+1})$.
Therefore, given $q_0$ and $q_1$ (or $q(t_0)$ and $\dot q(t_0)$) equations (\ref{back:pm1})-(\ref{back:pm2}) can be solved iteratively to find $q_2,\dots,q_n$.
Note that propagating the system in this manner ensures that the variational relation described in equation (\ref{back:actionsum}) is satisfied. 

The variational integrator is implemented through the introduction of the integration equation 
\begin{align}\label{back:fint}
f(q_{k+1}) = p_k + D_1L_\textrm{d}(q_k,q_{k+1}) = 0. 
\end{align}   
The Newton--Raphson method, outlined in Algorithm \ref{algo:1}, is used to find an approximate solution of equation (\ref{back:fint}).
Given $q_0$ and $q_1$ the integration scheme is initialized as 
\begin{align}\label{initalcond1}
p_1 =  D_2L_\textrm{d}(q_{0},q_{1}).
\end{align}  
Alternatively, given $q(t_0)$ and $\dot{q}(t_0)$ the integration scheme is initialized as 
\begin{align}\label{initalcond2}
p_1 = \frac{\partial}{\partial \dot{q}}L(q(t_0),\dot{q}(t_0)).
\end{align}  
The required derivatives can be found using the chain rule and equation (\ref{back:DL}) \cite{VI:Murphey:01}:
\begin{align}
D_1L_\textrm{d} &= \frac{h}{2} \frac{\partial L}{\partial q} - \frac{\partial L}{\partial \dot q}, \label{back:d1l}\\
D_2L_\textrm{d} &= \frac{h}{2} \frac{\partial L}{\partial q} + \frac{\partial L}{\partial \dot q},  \label{back:d2l}\\
Df(q_{k+1}) &= D_2D_1L_\textrm{d} = \frac{h}{4} \frac{\partial^2L}{\partial q \partial q}
+\frac{1}{2}\frac{\partial^2 L}{\partial \dot q\partial q}
-\frac{1}{2}\frac{\partial^2 L}{\partial q\partial \dot q}
-\frac{1}{h}\frac{\partial^2 L}{\partial \dot q\partial \dot q}.
\end{align}  

\begin{algorithm}[H]
\caption{Newton--Raphson Root Finder}
\label{algo:1}
$q_{k+1} = q_{k}$
\begin{algorithmic}
\While {$|f(q_{k+1})|>\epsilon_{\textrm{tol}}$}
    \State{$q_{k+1} \leftarrow q_{k+1}  -Df^{-1}(q_{k+1})\cdot f(q_{k+1})$}
\EndWhile
\end{algorithmic}
\end{algorithm}

\subsection{External Forces and Holonomic Constraints}
External forces can also be incorporated into the derivation of the variational integrator.
The Lagrange-d'Alembert principle is used to generalize the continuous Euler-Lagrange equation by modifying the variation of the action, $\delta S$, to
 \begin{align} \label{back:contELforce}
\delta S[q(t)] = \delta \int_{t_0}^{t_\textrm{f}}L(q(\tau), \dot q(\tau)) ~ \textrm{d} \tau + \int_{t_0}^{t_\textrm{f}}F(q(\tau), \dot q(\tau),u(\tau)) \cdot \delta q ~ \textrm{d} \tau
\end{align}   
where $F(q(\tau), \dot q(\tau),u(\tau))$ represents the total external forcing acting on the system and $u$ is the system's input (if any). 
Similar to the discretization of the Lagrangian, the left, $F_\textrm{d}^{-}(q_k,q_{k+1},u_k)$, and right, $F_\textrm{d}^+(q_k,q_{k+1},u_k)$, discrete forces are introduce in order to obtain a discrete equivalent to equation (\ref{back:contELforce}). 
The variation of the continuous external force is approximated over a small time interval as
 \begin{align} \label{dis:forces}
F_\textrm{d}^{-}(q_k,q_{k+1},u_k)\cdot\delta q_k + F_\textrm{d}^{+}(q_k,q_{k+1},u_k)\cdot\delta q_{k+1} \approx \int_{t_k}^{t_{k+1}}F(q(\tau), \dot q(\tau),u(\tau)) \cdot \delta q ~ \textrm{d} \tau
\end{align}     
where a midpoint approximation can be used to define the the left and right discrete forces as
 \begin{align}
F_\textrm{d}^{\pm}(q_k,q_{k+1},u_k)  &= \frac{h}{2}F( \frac{q_k +  q_{k+1}}{2}, \frac{q_{k+1}-q_k}{h},u_k),  
\end{align}   
and $u_k = u(t_k)$. 
The variational relation given in equation (\ref{back:actionsum}) can then be modified and the resulting forced DEL equation is given its \emph{position-momentum} form as   
\begin{align}
&p_{k} + D_1L_\textrm{d}(q_k,q_{k+1}) + F_\textrm{d}^{-}(q_k,q_{k+1},u_k)= 0, \label{back:pm1_f}\\
&p_{k} =  D_2L_\textrm{d}(q_{k-1},q_{k}) + F_\textrm{d}^{+}(q_{k-1},q_{k},u_{k-1}). \label{back:pm2_f}
\end{align}
The integrator equation and its derivative are now defined as
\begin{align}
f(q_{k+1}) &= p_{k} + D_1L_\textrm{d}(q_k,q_{k+1}) + F_\textrm{d}^{-}(q_k,q_{k+1},u_k), \label{back:forced:inteq}\\
Df(q_{k+1}) &= D_2D_1L_\textrm{d}(q_k,q_{k+1}) + D_2F_\textrm{d}^{-}(q_k,q_{k+1},u_k).
\end{align}
As before, given $q_0$, $q_1$, and the control input, $u(t)$, equation (\ref{back:forced:inteq}) can be solved iteratively to find $q_2,\dots,q_n$.
Given $q_0$ and $q_1$ the integration scheme is initialized as 
\begin{align}\label{initalcond1}
p_1 =  D_2L_\textrm{d}(q_{0},q_{1}) + F_\textrm{d}^{+}(q_{k-1},q_{k},u_{k-1}).
\end{align}  
Alternatively, given $q(t_0)$ and $\dot{q}(t_0)$ the integration scheme is initialized by equation (\ref{initalcond2}). 

Holonomic constraints can also be incorporated into the presented variational integrator. 
Specifically, the considered constraints are of the form $c(q) = [c_1(q), \dots, c_m(q)]^\textrm{T}$ where the system configuration is said to be valid if $c(q) = 0$.
Holomonic constraints restrict the set of possible system configurations to lie in a sub-manifold.
Therefore, during  propagation the computed system configurations should lie in the desired sub-manifold.
The integrator equation and its derivative can be modified to incorporate holonomic constraints \cite{marsden1999introduction}:
\begin{align}
f(q_{k+1},\lambda_k)&= \left[ \begin{array}{c} p_{k} + D_1L_\textrm{d}(q_k,q_{k+1}) + F_\textrm{d}^{-}(q_k,q_{k+1},u_k) - Dc^\textrm{T}(q_k)\lambda_k\\ c(q_{k+1}) \end{array} \right], \label{eqn:int:const1}\\
Df(q_{k+1},\lambda_k)&= \left[ \begin{array}{cc} D_2D_1L_\textrm{d}(q_k,q_{k+1}) + D_2F_\textrm{d}^{-}(q_k,q_{k+1},u_k)& - Dc^\textrm{T}(q_k) \\ Dc(q_{k+1})&0 \end{array} \right].\label{eqn:int:const2}
\end{align}
The term $Dc^\textrm{T}(q_k)\lambda_k$ represents a discretized force that imposes the constraint and $\lambda_k$ is the discrete Lagrange multiplier that defines the magnitude of this force. 
Note that the inclusion of the equation $c(q_{k+1}) = 0$ ensures that each discrete system configuration, $q_k$, observes the defined holomonic constraints.
The simple root finder algorithm in Algorithm \ref{algo:1} is modified such that the estimates of the discrete Lagrangian multipliers are also updated:
 \begin{align}
\left[ \begin{array}{c} q_{k+1}\\  \lambda_k \end{array} \right]  \leftarrow \left[ \begin{array}{c} q_{k+1}\\  \lambda_k \end{array} \right]  -Df^{-1}(q_{k+1},\lambda_k)\cdot f(q_{k+1},\lambda_k).
\end{align}

\section{Backward Error Analysis}
Consider an ordinary differential equation
 \begin{align}
\dot{x}(t) = f(x(t)), \quad x(0) = x_0, 
\end{align}
and the discrete propagation $x_\textrm{d}=\{x_0,x_1,\dots,x_{n}\}$ generated by a numerical method
 \begin{align}
x_{k+1} = \Psi(x_k),
\end{align}
that approximates the trajectory produced by the vector field such that $x_k \approx x(k h)$. 
It is assumed that the numerical method is \emph{consistent} and, therefore, $\lim_{h\rightarrow0}\frac{x_k-x(k h)}{h}=0$. 
Backward error analysis is used to generate \emph{a modified differential equation} of the form
 \begin{align}
f_{\textrm{mod}}(x(t)) &= f(x(t)) + h f_2(x(t))+ h^2 f_3(x(t)) + \dots, \label{mod:diffeq}\\
\dot{\tilde{x}}(t) &= f_{\textrm{mod}}(\tilde{x}(t)), \quad \tilde{x}(0) = x_0, 
\end{align}
such that $x_k = \tilde{x}(k h)$. 
That is, the discrete propagation generated by the numerical method exactly captures the evolution of the modified differential equation. 
Therefore, by comparing the original and modified differential equations the propagation distortion introduced by the integration method can be quantified directly in terms of the model instead of in terms of the states. 
We note that equation (\ref{mod:diffeq}) may diverge and, as a result, should be truncated when preforming rigorous analysis. 
In the context of this paper, convergence of (\ref{mod:diffeq}) is not of immediate concern and will not be examined closely. 

To derive the modified differential equation (up to a desired order) we first compute the Taylor series expansion of  $\tilde{x}(t + h)$ for a fixed $t$
\begin{align}
 \tilde{x}(t + h) &=  \tilde{x}(t)  + h(f(x) + h f_2(x(t)) + h^2 f_3(x(t)) + \dots) 
\nonumber \\&\quad+\frac{h^2}{2!}(f_{x}(x) + h f_{2,x}(x)+ \dots)(f(x) + h f_2(x(t)) + \dots)  
\nonumber \\&\quad+\frac{h^3}{3!}(f_{xx}(x) + h f_{2,xx}(x)+ \dots)\circ((f(x) + h f_2(x(t)) + \dots),(f(x) + h f_2(x(t)) + \dots)) + 
\nonumber \\&\quad +\frac{h^3}{3!}(f_{x}(x) + h f_{2,x}(x)+ \dots)(f_{x}(x) + h f_{2,x}(x)+ \dots)(f(x) + h f_2(x(t)) + \dots)\dots.
\end{align}
Next, it is assumed that the numerical method can be expanded as
 \begin{align}
\Psi(x) = x + h f(x) + h^2\psi_2(x) + h^3\psi_3(x) + \dots.
\end{align} 
Note that $f(x) = \psi_1(x)$ since the numerical method is consistent. 
Equating terms in the same power of $h$ gives the following recursive relations 
 \begin{align}
f_2(x) &= \psi_2(x) -\frac{1}{2!}f_{x}f, \\
f_3(x) &= \psi_3(x)-\frac{1}{3!}(f_{xx}\circ(f,f) + f_{x}f_{x}f) -\frac{1}{2!}(f_{x}f_2 + f_{2,x}f), \\
f_4(x) &= \dots .
\end{align}
Therefore, the accuracy of a numerical method can be directly quantified through the difference between $f(x(t))$ and $f_{\textrm{mod}}(x(t))$.
Generally, the order of the numerical method determines which modification terms are nonzero. 
For example, when analyzing a second order method it is expected that $f_2(x)=0$.  
References \cite{reich1999backward,hairer1994backward,WARMING1974159,griffiths1986scope} give a complete treatment of modified equations and backward error analysis. 

\subsection{Modified Lagrangians} \label{sec:backward:VI}
When backward error analysis is applied to variational or symplectic integrators the distortion introduced by the integration method can be described by \emph{modified Lagrangians and Hamiltonians}.
That is, the effect of the integration method can be described by a change in the system's Lagrangian or Hamiltonian. 
Furthermore, the modified dynamical systems described by the modified quantity are Hamiltonian systems \cite{hairer1994backward, vermeeren2015modified}.  
Therefore, the geometric properties associated with Hamiltonian systems are present in the  modified dynamical systems.
The analysis presented in this section closely follows the derivation of the modified Lagrangian presented in \cite{vermeeren2015modified}. 
We begin the analysis by considering the discrete Lagrangian (\ref{back:DL}),
\begin{eqnarray}\label{modlag:DL}
L_\textrm{d}(q_k,q_{k+1}) = h L\big(\frac{q_{k}+q_{k+1}}{2},\frac{q_{k+1}-q_{k}}{h}\big)\approx \int_{t_k}^{t_{k+1}}L(q(t),\dot{q}(t)) ~ \textrm{d} t.
\end{eqnarray}
Next, a Taylor series expansion around the midpoint, $q(\tau) =  q(\frac{t_{k+1}+t_k}{2})$, yields
\begin{eqnarray}
q_k &= q(\tau) -\frac{h}{2}\dot{q}(\tau)+\frac{h^2}{8}\ddot{q}(\tau) -\frac{h^3}{48}{q}^{(3)}(\tau) + o(h^4) \label{equ:qexp1} \\
q_{k+1} &= q(\tau) +\frac{h}{2}\dot{q}(\tau)+\frac{h^2}{8}\ddot{q}(\tau) +\frac{h^3}{48}{q}^{(3)}(\tau) + o(h^4) \label{equ:qexp2}
\end{eqnarray}
The discrete Lagrangian (\ref{modlag:DL}) can be equivalently given as a function of $q$ and its derivatives all evaluated at $\tau$\footnote{For ease of exposition, we denote $L( q(\tau) , \dot q(\tau) )$ as $L(\tau)$.},  
\begin{align}
\mathcal{L}\big(\tau\big) &= L( q(\tau) + \frac{h^2}{8}\ddot{q}(\tau) + o(h^4),\dot{q}(\tau) + \frac{h^2}{24}{q}^{(3)}(\tau)+ o(h^4)),\nonumber \\
&= L(\tau) + \frac{\partial L(\tau)}{\partial q}\Big( \frac{h^2}{8}\ddot{q}(\tau)+ o(h^4)\Big) +  \frac{\partial L(\tau)}{\partial \dot q}\Big(\frac{h^2}{24}{q}^{(3)}(\tau)+ o(h^4)\Big),\nonumber\\
  &= L(\tau) + \frac{h^2}{8}\frac{\partial L(\tau)}{\partial q}\ddot{q}(\tau) +  \frac{h^2}{24}\frac{\partial L(\tau)}{\partial \dot q}{q}^{(3)}(\tau) + o(h^4). \label{equ:errorindl}
\end{align}
Note if a generalized midpoint approximation is used equation (\ref{equ:errorindl}) would contain terms of order $h$ and $h^3$. Lemma \ref{lemma} is now given to find the approximation of the integral of the discrete Lagrangian. 

\begin{lemma} [\cite{vermeeren2015modified}, Lemma 6]\label{lemma}
For a smooth function $T : \mathbb{R}\rightarrow\mathbb{R}^N$ the following holds, 
\begin{align}
 h T(\tau)= \int^{t_{k+1}}_{t_k} T(\tau)  \textrm{d} t =\int^{t_{k+1}}_{t_k} T(t) + \sum^\infty_{i=1} h^{2i}(2^{1-2i}-1)\frac{B_{2_{i}}}{(2i)!}T^{2i}(t) \textrm{d} t, \label{eqn:int_exp}
\end{align}
where $B_i$ are the Bernoulli numbers, $\tau = \frac{t_{k+1}+t_k}{2}$ and $h = {t_{k+1}- t_k}$.
\end{lemma}
The proof of Lemma \ref{lemma} depends on a relatively straightforward application of the Euler--Maclaurin formula. 

Now suppose that there exists a modified Lagrangian, $L_\textrm{m}$, such that the discrete $L_\textrm{d}(q_k,q_{k+1})$ equals $\int^{t_{k+1}}_{t_k} L_\textrm{m} \textrm{d} t$. 
That is, by approximating the action integral of the dynamical system described by the Lagrangian, $L$, $L_\textrm{d}(q_k,q_{k+1})$ captures, up to some order of accuracy, the action integral of another system described by a  modified Lagrangian, $L_\textrm{m}$.
The relationship between the the Lagrangian, $L$, and the modified Lagrangian, $L_\textrm{m}$, is obtained from equations (\ref{equ:errorindl}) and (\ref{eqn:int_exp}):
\begin{align}
\int^{t_{k+1}}_{t_k} L_\textrm{m} \ \textrm{d} t &= L_\textrm{d}(q_k,q_{k+1})= \int^{t_{k+1}}_{t_k} {L}\big(\frac{q_{k}+q_{k+1}}{2},\frac{q_{k+1}-q_{k}}{h}\big) \ \textrm{d} t= \int^{t_{k+1}}_{t_k} \mathcal{L}\big(\tau\big) \ \textrm{d} t \nonumber \\
&=\int^{t_{k+1}}_{t_k} L(t) + \frac{h^2}{8}\frac{\partial L(t)}{\partial q}\ddot{q}(t) +  \frac{h^2}{24}\frac{\partial L(t)}{\partial \dot q}{q}^{(3)}(t) - \frac{h^2}{24}\ddot{L}(t) + o(h^4) \ \textrm{d} t. \label{modequaltoexp}
\end{align}
An expression for the modified Lagrangian is obtained, 
\begin{align}
L_\textrm{m} &= L+ \frac{h^2}{8}\frac{\partial L}{\partial q}\ddot{q} +  \frac{h^2}{24}\frac{\partial {L}}{\partial \dot q}{q}^{(3)}  - \frac{h^2}{24}\ddot{L}(\tau) + o(h^4), 
\intertext{or, equivalently,}
L_\textrm{m} &= L - \frac{h^2}{24}(-2\frac{\partial L}{\partial q}\ddot{q} + \dot{q}^\textrm{T}\frac{\partial^2 L}{\partial q\partial q}\dot{q}+ \ddot{q}^\textrm{T}\frac{\partial^2 L}{\partial \dot q\partial \dot q}\ddot{q}+ 2\ddot{q}^\textrm{T}\frac{\partial^2 L}{\partial q\partial \dot q}\dot{q}) + o(h^4).  
\end{align}
Note that since the modified Lagrangian and the Lagrangian differ in terms of order $h^2$ the evolution of the modified system can be expressed as
\begin{eqnarray}\label{equ:modsys}
\ddot{q} = \Big(\frac{\partial^2 L_\textrm{m}}{\partial \dot q\partial \dot q}\Big)^{-1}\Big(\frac{\partial L_\textrm{m}}{\partial q} - \frac{\partial^2 L_\textrm{m}}{\partial q\partial \dot q}\dot q\Big)=\Big(\frac{\partial^2 L}{\partial \dot q\partial \dot q}\Big)^{-1}\Big(\frac{\partial L}{\partial q} - \frac{\partial^2 L}{\partial q\partial \dot q}\dot q\Big)+ o(h^2).
\end{eqnarray}
Equation (\ref{equ:modsys}) is used to express the modified Lagrangian as a function of $q$ and $\dot{q}$:
\begin{eqnarray} \label{equ:modL}
L_\textrm{m} =  {L} -\frac{h^2}{24}\Big(\dot{q}^T\big(\frac{\partial^2 {L}}{\partial q\partial q} - \frac{\partial^2 {L}}{\partial\dot q\partial q}^\textrm{T}\frac{\partial^2 {L}}{\partial\dot q\partial\dot q}^{-1}\frac{\partial^2 {L}}{\partial\dot q\partial q}\big)\dot{q} - \frac{\partial {L}}{q}^\textrm{T}\frac{\partial^2 {L}}{\partial\dot q\partial\dot q}^\textrm{-1}\frac{\partial {L}}{q}
+ 2\dot{q}^\textrm{T}\frac{\partial^2 {L}}{\partial\dot q\partial q}^\textrm{T}\frac{\partial^2{L}}{\partial\dot q\partial\dot q}^{-1}\frac{\partial {L}}{\partial q}\Big) + o(h^4).
\end{eqnarray}
The modified Lagrangian in equation (\ref{equ:modL}) quantifies the propagation distortion of order $h^2$ introduced by the midpoint variational integrator. 
Note that the induced distortion is a function of the partial derivatives of the Lagrangian. 
As a result, the modification of any sufficiently differentiable Lagrangian can be easily computed. 
Section \ref{sec:higher} extends the presented analysis in order to quantify modifications for any desired order.  

\section{Surrogate Lagrangians}
In the previous section backward error analysis quantified the error induced by a numerical integration scheme by defining a modified differential equations.
Intuitively, if the manner in which a integration method ``modifies" a differential equation is known then a new integration scheme can be constructed that mitigates (or eliminates) known errors. 
This idea has led to the creation of integration methods for general ordinary differential equations and special classes of Hamiltonian systems \cite{mushtaq2014higher,Mushtaq20141461,chartier2007numerical,hairer2006preprocessed,StochModEqu,chartier2007modified,kozlov2008high}.
In this section, \emph{surrogate Lagrangians} are derived and analyzed.  
It is shown that the accuracy of a midpoint variational integrator is improved if the considered system's Lagrangian is replaced by its surrogate Lagrangian. 
To begin, define the second-order \emph{surrogate} Lagrangian, $\hat{L}$, as
\begin{eqnarray} \label{surr:lag}
\hat{L} =  {L} + \frac{h^2}{24}(-2\frac{\partial L}{\partial q}\ddot{q} + \dot{q}^\textrm{T}\frac{\partial^2 L}{\partial q\partial q}\dot{q}+ \ddot{q}^\textrm{T}\frac{\partial^2 L}{\partial \dot q\partial \dot q}\ddot{q}+ 2\ddot{q}^\textrm{T}\frac{\partial^2 L}{\partial q\partial \dot q}\dot{q}).
\end{eqnarray}
Equation (\ref{equ:modL}) is used to derive an expression for the the\emph{ modified surrogate} Lagrangian, $\hat{L}_\textrm{m}$:
\begin{align} 
\hat{L}_\textrm{m} &=  \hat{L} - \frac{h^2}{24}(-2\frac{\partial \hat{L}}{\partial q}\ddot{q} + \dot{q}^\textrm{T}\frac{\partial^2 \hat{L}}{\partial q\partial q}\dot{q}+ \ddot{q}^\textrm{T}\frac{\partial^2 \hat{L}}{\partial \dot q\partial \dot q}\ddot{q}+ 2\ddot{q}^\textrm{T}\frac{\partial^2 \hat{L}}{\partial q\partial \dot q}\dot{q}).
\end{align}
The modified surrogate Largrangian and the Lagrangian are related as,
\begin{align}
\hat{L}_\textrm{m}  &=  {L} + o(h^4), \label{equ:mod2}
\end{align}
by noting that
\begin{align} 
-2\frac{\partial \hat{L}}{\partial q}\ddot{q} + \dot{q}^\textrm{T}\frac{\partial^2 \hat{L}}{\partial q\partial q}\dot{q}+ \ddot{q}^\textrm{T}\frac{\partial^2 \hat{L}}{\partial \dot q\partial \dot q}\ddot{q}+ 2\ddot{q}^\textrm{T}\frac{\partial^2 \hat{L}}{\partial q\partial \dot q}\dot{q} = -2\frac{\partial L}{\partial q}\ddot{q} + \dot{q}^\textrm{T}\frac{\partial^2 L}{\partial q\partial q}\dot{q}+ \ddot{q}^\textrm{T}\frac{\partial^2 L}{\partial \dot q\partial \dot q}\ddot{q}+ 2\ddot{q}^\textrm{T}\frac{\partial^2 L}{\partial q\partial \dot q}\dot{q} + o(h^2).
\end{align}
Since the modified surrogate Lagrangian and the Lagrangian differ in terms of order $h^4$ the evolution of the modified surrogate system can be expressed as
\begin{eqnarray}\label{equ:modsurrsys}
\ddot{q} = \Big(\frac{\partial^2 \hat{L}_\textrm{m}}{\partial \dot q\partial \dot q}\Big)^{-1}\Big(\frac{\partial \hat{L}_\textrm{m}}{\partial q} - \frac{\partial^2 \hat{L}_\textrm{m}}{\partial q\partial \dot q}\dot q\Big)=\Big(\frac{\partial^2 L}{\partial \dot q\partial \dot q}\Big)^{-1}\Big(\frac{\partial L}{\partial q} - \frac{\partial^2 L}{\partial q\partial \dot q}\dot q\Big)+ o(h^4).
\end{eqnarray}
The modified surrogate differential equation approximates the original differential equation up to order $h^4$ while the expression given in equation (\ref{equ:modL}) only does so up to order $h^2$. 
The increase of accuracy was not achieved through an increase in the complexity of the variational integrator. 
In some sense, the modification to the Hamiltonian system made by the variational integrator was ``corrected" through the introduction of the surrogate system. 
Therefore, as shown in numerical examples in the following sections, a second order integrator can achieve fourth order accuracy.
Furthermore, as shown in Section \ref{sec:higher} the methodology can be extended in order to increase the integrator accuracy up to an arbitrary order. 

The surrogate variational integrator can now be defined by the DEL equations
\begin{align}
&p_{k} + D_1\hat{L}_\textrm{d}(q_k,q_{k+1}) = 0,\\
&p_{k} =  D_2\hat{L}_\textrm{d}(q_{k-1},q_{k}),
\end{align}
and the initial conditions given by
\begin{align}
p_1 =  D_2L_\textrm{d}(q_{0},q_{1}) \quad \textrm{or} \quad p_1 = \frac{\partial}{\partial \dot{q}}L(q(t_0),\dot{q}(t_0)).
\end{align}  
Note that the initial conditions are defined with the nominal Lagrangian and not the surrogate Lagrangian. 

\subsection{Example: Harmonic Oscillator}\label{sec:lowharmosc}
\begin{figure}[h]
    \centering
    \begin{subfigure}[]{0.5\columnwidth}
	\centering
        	\includegraphics[width=\columnwidth]{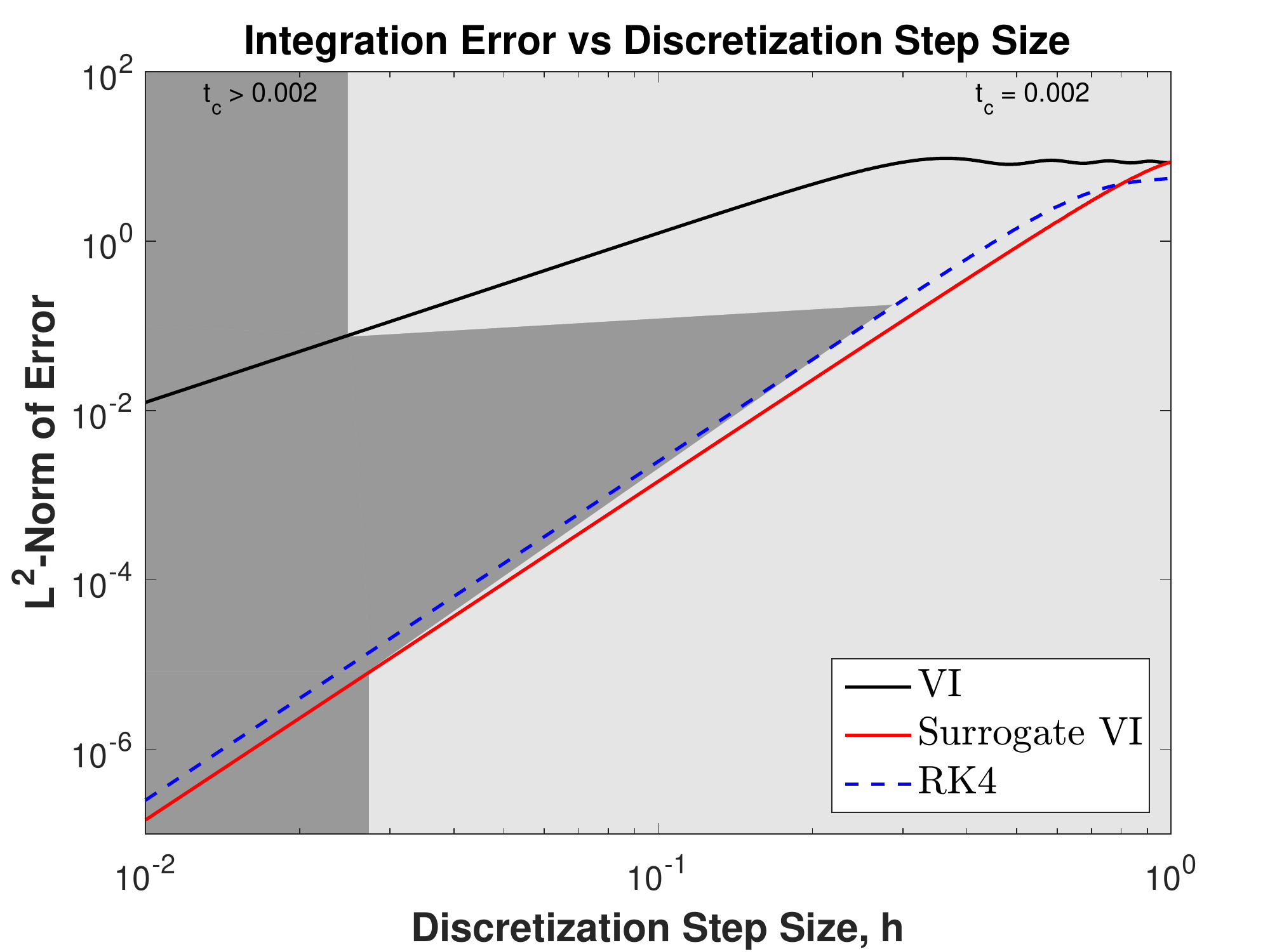}
	\caption{ }
    \end{subfigure}%
    \begin{subfigure}[]{0.5\columnwidth}
	\centering
        	\includegraphics[width=\columnwidth]{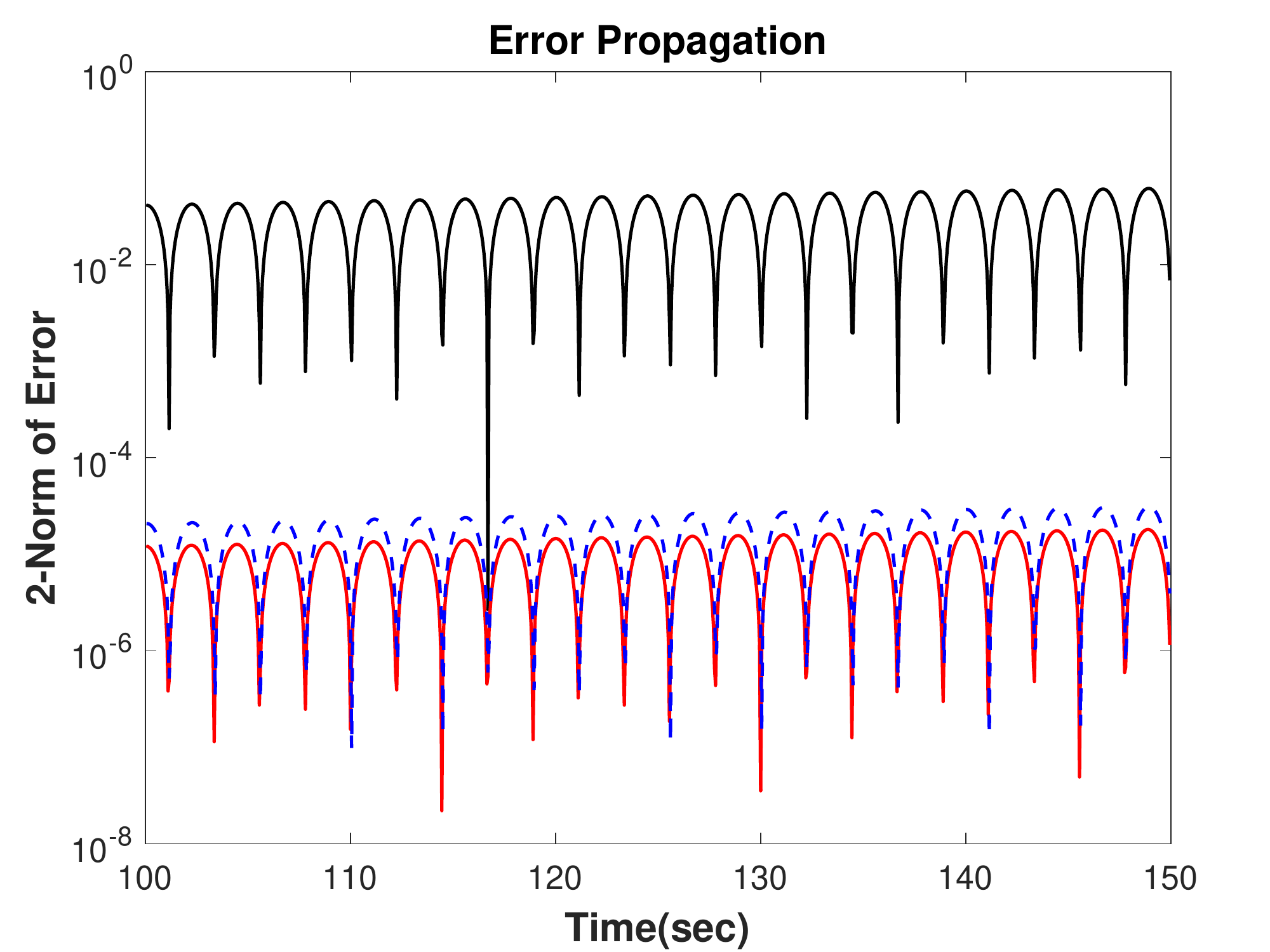}
	\caption{ }
    \end{subfigure}
    \caption{(a): The $L^2$-norm of the error as a function of the utilized discretization time step plotted on a contour of average computational times. 
The Runge--Kutta fourth-order method and the surrogate variational integrator exhibit a fourth order convergence of the $L^2$-norm of the error.
The increase level of accuracy in the surrogate variational integrator is not associated with a greater level of computational effort.
(b): The 2-norm of error for a particular execution of the numerical integrators when $h=0.05$ in the time interval $100\leq t \leq 150$.
The surrogate variational integrator tends to be the most accurate numerical method. }
    \label{fig:spring}
\end{figure}
In order to further illustrate our discussion of surrogate Lagrangians and variational integrators, consider a mass-spring system with mass $M$ and spring constant $K$.
Its Lagrangian is given as 
\begin{eqnarray}
{L}&=  \frac{1}{2}M\dot q^2 - \frac{1}{2}K q^2,
\end{eqnarray}
and its surrogate Lagrangian is given as
\begin{eqnarray}
\hat{L}&=  \frac{1}{2}(M-K\frac{h^2}{12})\dot q^2 - \frac{1}{2}(K+K{M}^{-1}K\frac{h^2}{12})q^2.
\end{eqnarray}
Note that the surrogate Lagrangian simply describes another mass-spring system with a different mass and spring constant. 
Furthermore, the surrogate system's mass and spring constant are dependent on the discretization time step and $\hat{L}\rightarrow{L}$ as $h\rightarrow 0$. 
The system was propagated using a variational integrator, a surrogate variational integrator, and the classical fourth order Runge-–Kutta method. 
The initial condition of the system was set as $q(0)=0$ and $\dot{q}(0)=1$ and $M=1$ and $K=2$. 
The system was propagated for 150 seconds. 
We define the 2-norm of error, $e_2(t)$, and the $L^2$-norm of the error, $e_{L^2}$, as 
\begin{align} 
e_2(t)  &= \sqrt{(q(t) - q_\textrm{a}(t))^2}, \\
e_{L^2}  &= \Big(\int^{t_\textrm{f}}_0(q(t) - q_\textrm{a}(t))^2  \ \textrm{d} t\Big)^{\frac{1}{2}}
\end{align}
where $q_\textrm{a}(t)$ is the trajectory obtained from the known analytic solution. 
Figure \ref{fig:spring} shows the 2-norm error when $h=0.05$ and $L^2$-norm of the error for a range of discretization time steps.
A contour of the average computational time from 100 executions of each integrator configuration is also displayed.
Note that both the Runge--Kutta fourth order method and the surrogate variational integrator exhibit a fourth order convergence of the $L^2$-norm of the error. 
The surrogate variational integrator is more accurate than the Runge--Kutta method for most of the discretization time steps considered. 
Furthermore, the surrogate variational integrator is able to achieve the lowest $L^2$-norm in each contour depicted. 
Therefore, the improvement of accuracy is not accompanied, in this case, by an increase in computational effort.

\subsection{Holonomic Constraints}
In this section, the class of Hamiltonian systems in which a surrogate Lagrangian can be obtained is expanded by considering holonomic constraints.
The analysis presented in this section outlines arguments originally used to obtained a modified Lagrangian for constrained systems presented in \cite{vermeeren2015modified}.
To begin, the augmented Lagrangian 
\begin{align}
\bar{L}(q, \dot{q}) = {L}(q, \dot{q}) + \lambda^\textrm{T}c(q),
\end{align}
is introduced to derive the Euler-Lagrange equation of the constrained system:
\begin{align} 
\ddot{q} =\Big(\frac{\partial^2\bar L}{\partial \dot q\partial \dot q}\Big)^{-1}\Big(\frac{\partial \bar L}{\partial q} - \frac{\partial^2 \bar L}{\partial \dot q\partial q}\dot q\Big)
 =\Big(\frac{\partial^2 L}{\partial \dot q\partial \dot q}\Big)^{-1}\Big(\frac{\partial  L}{\partial q} +\frac{\partial c}{\partial q}^\textrm{T} \lambda - \frac{\partial^2  L}{\partial \dot q\partial q}\dot q\Big).
\end{align}
The Lagrange multiplier $\lambda$ can be replaced by an explicit function of $q$ and $\dot{q}$. 
This explicit function is found by taking the time derivative of $c(q) = 0$ twice and then solving for $\lambda$. 
As a result,  the Euler-Lagrange equations can be rewritten as a function of $q$ and $\dot{q}$: 
\begin{align} \label{const:syseq}
\ddot{q} =\Big(\frac{\partial^2 L}{\partial \dot q\partial \dot q}\Big)^{-1}\Big(\frac{\partial L}{\partial q} + \frac{\partial c}{\partial q}^\textrm{T}\lambda(q, \dot{q}) - \frac{\partial^2 L}{\partial \dot q\partial q}\dot q\Big).
\end{align}
An expression for the surrogate Lagrangian is derived from equations (\ref{surr:lag}) and (\ref{const:syseq}):
\begin{align} 
\hat{L}&=  {L} +\frac{h^2}{24}\Big(\dot{q}^T\big(\frac{\partial^2 {L}}{\partial q\partial q} - \frac{\partial^2 {L}}{\partial\dot q\partial q}^\textrm{T}\frac{\partial^2 {L}}{\partial\dot q\partial\dot q}^{-1}\frac{\partial^2 {L}}{\partial\dot q\partial q}\big)\dot{q} - \frac{\partial {L}}{q}^\textrm{T}\frac{\partial^2 {L}}{\partial\dot q\partial\dot q}^\textrm{-1}\frac{\partial {L}}{q}
\nonumber\\ &\quad+ 2\dot{q}^\textrm{T}\frac{\partial^2 {L}}{\partial\dot q\partial q}^\textrm{T}\frac{\partial^2{L}}{\partial\dot q\partial\dot q}^{-1}\frac{\partial {L}}{\partial q} + \lambda(q, \dot{q})^\textrm{T} \frac{\partial c}{\partial q}\frac{\partial^2 {L}}{\partial\dot q\partial\dot q}^{-1} \frac{\partial c}{\partial q}^\textrm{T}\lambda(q, \dot{q}) \Big) + o(h^4).
\end{align}
The surrogate system should evolve in the same constrained manifold as the original system. 
Therefore, the considered holonomic constraint is not changed and the augmented surrogate Lagrangian is given as
\begin{align}
\bar{\hat{L}}(q, \dot{q}) = \hat{L}(q, \dot{q}) + \hat\lambda^\textrm{T}c(q)
\end{align}
where $\hat\lambda$ is the surrogate Lagrange multiplier.

\subsection{External Forces}
External forces change the expression of the surrogate Lagrangian in a similar manner as holonomic constraints. 
As before, the forced Euler-Lagrange equation (\ref{back:contELforce}) and equation (\ref{surr:lag}) are used to derive an expression for the surrogate Lagrangian for the forced case: 
\begin{align} 
\hat{L}&=  {L} +\frac{h^2}{24}\Big(\dot{q}^T\big(\frac{\partial^2 {L}}{\partial q\partial q} - \frac{\partial^2 {L}}{\partial\dot q\partial q}^\textrm{T}\frac{\partial^2 {L}}{\partial\dot q\partial\dot q}^{-1}\frac{\partial^2 {L}}{\partial\dot q\partial q}\big)\dot{q} - \frac{\partial {L}}{q}^\textrm{T}\frac{\partial^2 {L}}{\partial\dot q\partial\dot q}^\textrm{-1}\frac{\partial {L}}{q}
\nonumber\\ &\quad+ 2\dot{q}^\textrm{T}\frac{\partial^2 {L}}{\partial\dot q\partial q}^\textrm{T}\frac{\partial^2{L}}{\partial\dot q\partial\dot q}^{-1}\frac{\partial {L}}{\partial q} + F(q,\dot q,u)^\textrm{T}\frac{\partial^2 {L}}{\partial\dot q\partial\dot q}^{-1} F(q,\dot q,u) \Big) + o(h^4).
\end{align}
However, unlike holonomic constraints, external forces are approximated by the variational integrator using left, $F_\textrm{d}^{-}(q_k,q_{k+1},u_k)$, and right, $F_\textrm{d}^+(q_k,q_{k+1},u_k)$, discrete forces. 
Therefore, in this section backward error analysis is used to find the modified external force.
That is, the modified external force and the modified Lagrangian characterize a forced system that is exactly captured by the implemented variational integrator. 
An expression for the surrogate external force is then found. 
To begin, consider the Taylor series expansion around the midpoint, $q(\tau) =  q(\frac{t_{k+1}+t_k}{2})$, given by equations (\ref{equ:qexp1}) and (\ref{equ:qexp2}).
Furthermore, the system's input, $u(t)$, is expanded as 
\begin{eqnarray}
u_k &= u(\tau) -\frac{h}{2}\dot{u}(\tau)+\frac{h^2}{8}\ddot{u}(\tau) -\frac{h^3}{48}{u}^{(3)}(\tau) + o(h^5),\nonumber \\
u_{k+1} &= u(\tau) +\frac{h}{2}\dot{u}(\tau)+\frac{h^2}{8}\ddot{u}(\tau) +\frac{h^3}{48}{u}^{(3)}(\tau) + o(h^5).\nonumber
\end{eqnarray}
As before, the discrete forces (\ref{dis:forces}) can be equivalently given as functions of $q$, $u$, and their derivatives all evaluated at $\tau$
 \footnote{The definition of the discrete forces were changed from $\frac{h}{2}F(\frac{q_k + q_{k+1}}{2}, \frac{q_{k+1}-q_k}{h},u_k)$ to $\frac{h}{2}F(\frac{q_k + q_{k+1}}{2}, \frac{q_{k+1}-q_k}{h},\frac{u_k + u_{k+1}}{2})$ to avoid any terms of odd order ($h, h^3, \dots$). However, the analysis presented here can be done considering the original definition.},
\begin{align}
\mathcal{F}(\tau)\cdot\delta q_{k}(\tau) =& \frac{1}{2}F(q(\tau) + \frac{h^2}{8}\ddot{q}(\tau) \dots,\dot{q}(\tau) + \frac{h^2}{24}{q}^{(3)}(\tau)\dots, u(\tau) + \frac{h^2}{8}\ddot{u}(\tau)\dots)\cdot(\delta q(\tau) -\frac{h}{2}\dot{\delta  q}(\tau)\dots)\nonumber \\ +
&\frac{1}{2}F(q(\tau) + \frac{h^2}{8}\ddot{q}(\tau) \dots,\dot{q}(\tau) + \frac{h^2}{24}{q}^{(3)}(\tau)\dots, u(\tau) + \frac{h^2}{8}\ddot{u}(\tau)\dots)\cdot(\delta q(\tau) +\frac{h}{2}\dot{\delta  q}(\tau)\dots)\nonumber \\ 
=& F(q(\tau) + \frac{h^2}{8}\ddot{q}(\tau) \dots,\dot{q}(\tau) + \frac{h^2}{24}{q}^{(3)}(\tau)\dots, u(\tau) + \frac{h^2}{8}\ddot{u}(\tau)\dots)\cdot(\delta q(\tau) + \frac{h^2}{24}{\delta q}^{(3)}(\tau)\dots)\nonumber \\ 
\mathcal{G}(\tau)=& G(\tau) + \frac{h^2}{8}\frac{\partial G(\tau)}{\partial \delta q} \delta \ddot{q}  + \frac{h^2}{8}\frac{\partial G(\tau)}{\partial q}\ddot{q}+ \frac{h^2}{8}\frac{\partial G(\tau)}{\partial u}\ddot{u}+ \frac{h^2}{24} \frac{\partial G(\tau)}{\partial \dot q}{q}^{(3)} + o(h^4) 
\end{align}
where $G(t)  = F(t) \cdot \delta q(t)$ is the virtual work done by force $F(t)$. 
As with the discrete Lagrangian, assume that there exists a modified force that is captured perfectly by the left and right discrete forces, 
\begin{align} 
\int_{t_k}^{t_{k+1}}F_\textrm{m}(t) \cdot \delta q(t) ~ \textrm{d} t = \frac{h}{2}(F^{-}_\textrm{d}\cdot\delta q_k + F^{+}_\textrm{d}\cdot\delta q_{k+1})=  \int_{t_k}^{t_{k+1}}\mathcal{F}(\tau)\cdot\delta q_{k}(\tau)~ \textrm{d} t.
\end{align}
Lemma \ref{lemma} is used to derive the expression for the modified virtual work, $G_\textrm{m}$, as
\begin{align}
G_\textrm{m} =& G(\tau) + \frac{h^2}{8}\frac{\partial G(\tau)}{\partial \delta q} \delta \ddot{q}  + \frac{h^2}{8}\frac{\partial G(\tau)}{\partial q}\ddot{q}+ \frac{h^2}{8}\frac{\partial G(\tau)}{\partial u}\ddot{u}+ \frac{h^2}{24} \frac{\partial G(\tau)}{\partial \dot q}{q}^{(3)} - \frac{h^2}{24}\ddot{G}(\tau) + o(h^4) \nonumber\\
=& G - \frac{h^2}{24}(
-2\frac{\partial G}{\partial q}\ddot{q} 
-2\frac{\partial G}{\partial \delta q}\delta\ddot{ q} 
-2\frac{\partial G}{\partial \delta u}\delta\ddot{u} 
+\dot{q}^\textrm{T}\frac{\partial^2 G}{\partial q\partial q}\dot{q}
+\dot{u}^\textrm{T}\frac{\partial^2 G}{\partial u\partial u}\dot{u}
+\ddot{q}^\textrm{T}\frac{\partial^2 G}{\partial \dot q\partial \dot q}\ddot{q}
+2\delta\dot{q}^\textrm{T}\frac{\partial^2 G}{\partial q\partial \delta   q}\dot{q} \nonumber\\
&+2\dot{u}^\textrm{T}\frac{\partial^2 G}{\partial q\partial u}\dot{q}
+2\ddot{q}^\textrm{T}\frac{\partial^2 G}{\partial q\partial \dot q}\dot{q}
+2\dot{u}^\textrm{T}\frac{\partial^2 G}{\partial \delta q\partial u}\delta\dot{q}
+2\ddot{q}^\textrm{T}\frac{\partial^2 G}{\partial \delta q\partial \dot q}\delta\dot{q}
+2\ddot{q}^\textrm{T}\frac{\partial^2 G}{\partial u \partial \dot q}\dot{u}
)+ o(h^4).  
\end{align}
Note that $\frac{\partial^2 G}{\partial \delta q\partial \delta q} = 0$. 
Since the modified Lagrangian and modified virtual work only differ in terms of order $h^2$ from the modeled ones the evolution of the modified system and the associated virtual displacement can be expressed as
\begin{align}
\ddot{q} =& \Theta(q,\dot{q},u) + o(h^2), \label{equ:fqom} \\ 
\delta\ddot{q} =& D_1\Theta(q,\dot{q},u)\delta q+ D_2 \Theta(q,\dot{q},u)\delta \dot q+ o(h^2),\label{equ:dfqom}
\end{align}
where $\Theta$ is the equations of motions obtained from the Lagrange-d'Alembert principle. 
Equations (\ref{equ:fqom}) and (\ref{equ:dfqom}) are used to eliminate any dependence of $G_\textrm{m}$ on $\ddot q$ and $\delta \ddot q$.
$G_\textrm{m}$ is rewritten in the following form
\begin{align}
G_\textrm{m} = H_1(q,\dot q, u)\cdot \delta q + H_2(q,\dot q, u)\cdot \delta\dot  q + o(h^4). 
\end{align}
Integration by parts yields an expression for the modified force: 
\begin{align}\label{int:byparts}
\int_{t_{k}}^{t_{k+1}}G_\textrm{m}~\textrm{d} t &= H_2\cdot \delta q|^{t_{k+1}}_{t_k} +  \int_{t_{k}}^{t_{k+1}}(H_1 - \dot{H}_2 + o(h^4))\cdot \delta q ~\textrm{d} t.
\end{align}
The surrogate virtual work is now derived as
\begin{align}
\hat{G} &= G + \frac{h^2}{24}(
-2\frac{\partial G}{\partial q}\ddot{q} 
-2\frac{\partial G}{\partial \delta q}\delta\ddot{ q} 
-2\frac{\partial G}{\partial u}\delta\ddot{u} 
+\dot{q}^\textrm{T}\frac{\partial^2 G}{\partial q\partial q}\dot{q}
+\dot{u}^\textrm{T}\frac{\partial^2 G}{\partial u\partial u}\dot{u}
+\ddot{q}^\textrm{T}\frac{\partial^2 G}{\partial \dot q\partial \dot q}\ddot{q}
+2\delta\dot{q}^\textrm{T}\frac{\partial^2 G}{\partial q\partial \delta   q}\dot{q} \nonumber\\
&+2\dot{u}^\textrm{T}\frac{\partial^2 G}{\partial q\partial u}\dot{q}
+2\ddot{q}^\textrm{T}\frac{\partial^2 G}{\partial q\partial \dot q}\dot{q}
+2\dot{u}^\textrm{T}\frac{\partial^2 G}{\partial \delta q\partial u}\delta\dot{q}
+2\ddot{q}^\textrm{T}\frac{\partial^2 G}{\partial \delta q\partial \dot q}\delta\dot{q}
+2\ddot{q}^\textrm{T}\frac{\partial^2 G}{\partial u \partial \dot q}\dot{u}
)+ o(h^4),
\end{align}
and can be equivalently represented as 
\begin{align}
\hat{G} &= \hat{H}_1(q,\dot q, u)\cdot \delta q + \hat{H}_2(q,\dot q, u)\cdot \delta\dot  q.
\end{align}
It is easy to verified that 
\begin{align}
\hat{G}_\textrm{m}  = G + o(h^4).
\end{align}
Referring to equation (\ref{int:byparts}), the surrogate virtual work is then discretized as   
\begin{align}
\int_{t_{k}}^{t_{k+1}}\hat{G}_\textrm{}~\textrm{d} t =& \hat{H}_2(q_{k+1}, \frac{q_{k+1}-q_k}{h}, u_{k+1})\cdot \delta q({t_{k+1}}) -\hat{H}_2(q_k, \frac{q_{k}-q_{k-1}}{h},u_k)\cdot \delta q({t_{k}}) \nonumber\\&
+  (\hat{H}^{+}_{1,\textrm{d}} - \dot{\hat{H}}^{+}_{2,\textrm{d}})\cdot \delta q({t_{k+1}}) +   (\hat{H}^{-}_{1,\textrm{d}} - \dot{\hat{H}}^{-}_{2,\textrm{d}})\cdot \delta q({t_{k}}),
\end{align}
where 
\begin{align}
\hat{H}^{\pm}_{i,\textrm{d}} = \frac{h}{2}\hat{H}(\frac{q_{k+1} + q_{k}}{2}, \frac{q_{k+1} - q_{k}}{h}, \frac{u_{k+1} + u_{k}}{2}).  
\end{align}
As a result, the surrogate discrete force is expressed as
\begin{align}\label{equ:surroodiscreteforce}
\hat{F}_\textrm{d}^{-}(q_{k-1},q_k,q_{k+1},u_k,u_{k+1}) &= -\hat{H}_2(q_k, \frac{q_{k}-q_{k-1}}{h},u_k) +  \hat{H}^{-}_{1,\textrm{d}} - \dot{\hat{H}}^{-}_{2,\textrm{d}},\nonumber \\
\hat{F}_\textrm{d}^{+}(q_k,q_{k+1},u_k,u_{k+1}) &= \hat{H}_2(q_{k+1}, \frac{q_{k+1}-q_k}{h}, u_{k+1})+  \hat{H}^{+}_{1,\textrm{d}} - \dot{\hat{H}}^{+}_{2,\textrm{d}}.
\end{align}
The resulting forced \emph{surrogate} DEL equations are now given as   
\begin{align}
&p_{k} + D_1\hat{L}_\textrm{d}(q_k,q_{k+1}) + \hat{F}_\textrm{d}^{-}(q_{k-1}, q_k,q_{k+1},u_k,u_{k+1})= 0, \label{back:pm1new}\\
&p_{k} =  D_2\hat{L}_\textrm{d}(q_{k-1},q_{k}) + \hat{F}_\textrm{d}^{+}(q_{k-1},q_{k},u_{k-1}, u_k). \label{back:pm2new}
\end{align}

\section{Higher Order Surrogate Lagrangians}\label{sec:higher}
Using a similar procedure as in Section \ref{sec:backward:VI}, expressions for modified Lagrangians that explicitly contain higher order terms can be derived.  
These expressions can then be used to define higher order surrogate Lagrangians.
Therefore, the modified differential equations resulting from these higher order surrogate Lagrangians approximate the modeled Euler-Lagrangian differential equations up to higher orders ($h^6, h^8,$ etc.).  
To begin, the Taylor series expansion around the midpoint, $q(\tau) =  q(\frac{t_{k+1}+t_k}{2})$, is extended, 
\begin{eqnarray}
q_k &=& q(\tau) -\frac{h}{2}\dot{q}(\tau)+\frac{h^2}{8}\ddot{q}(\tau) -\frac{h^3}{48}{q}^{(3)}(\tau) +\frac{h^4}{384}{q}^{(4)}(\tau) -\frac{h^5}{3840}{q}^{(5)}(\tau)+ o(h^6), \label{equ:highqexp1} \\
q_{k+1} &=& q(\tau) +\frac{h}{2}\dot{q}(\tau)+\frac{h^2}{8}\ddot{q}(\tau) +\frac{h^3}{48}{q}^{(3)}(\tau) +\frac{h^4}{384}{q}^{(4)}(\tau) +\frac{h^5}{3840}{q}^{(5)}(\tau)+ o(h^6). \label{equ:highqexp2}
\end{eqnarray}
This Taylor series expansion is used to obtain an expression for the discretization of the surrogate Lagrangian given in equation (\ref{surr:lag}):
\begin{align}
\hat{\mathcal{L}}(\tau)  =& \hat{L}\Big( q(\tau) + \frac{h^2}{8}\ddot{q}(\tau)+ \frac{h^4}{384}{q}^{(4)}(\tau) + o(h^6),\dot{q}(\tau) + \frac{h^2}{24}{q}^{(3)}(\tau)+ \frac{h^4}{1920}{q}^{(5)}(\tau) + o(h^6)\Big), \nonumber \\
 =& \hat{L}(\tau) + {h^2}\theta_2(\hat{L}(\tau))+ {h^4}\theta_4(\hat{L}(\tau))+ o(h^6),
\end{align}
where
\begin{align}
 \theta_2(L(\tau)) =& \frac{\partial L(\tau)}{\partial q}\frac{\ddot{q}}{8} + \frac{\partial L(\tau)}{\partial \dot q}\frac{{q}^{(3)}}{24}, \nonumber \\
 \theta_4(L(\tau)) =& \frac{\partial L(\tau)}{\partial q}\frac{{q}^{(4)}}{384}  + \frac{\partial L(\tau)}{\partial \dot q}\frac{{q}^{(5)}}{1920}+ \frac{\partial^2 L(\tau)}{\partial^2 q}\frac{\ddot{q}^2}{128}+ \frac{\partial^2 L(\tau)}{\partial^2 \dot q}\frac{({q}^{(3)})^2}{1152}+ \frac{\partial^2 {L}(\tau)}{\partial q\partial \dot q}\frac{{q}^{(3)} \ddot{q}}{96}.
\end{align}
As done in  Section \ref{sec:backward:VI}, Lemma \ref{lemma} is used to obtain an expression for the modified surrogate Lagrangian: 
\begin{align}
\hat{L}_\textrm{m}= \hat{L} + {h^2}\theta_2(\hat{L})+ {h^4}\theta_4(\hat{L}) - \frac{h^2}{24}(\ddot{\hat{L}} + h^2\ddot{\theta}_2(\hat{L})) + \frac{7h^4}{5760}\hat{L}^{(4)} + o(h^6).\label{equ:pre4mod}
\end{align}
Noting that $\hat{L} =L -h^2 \Phi_2(L)$ the modified surrogate Largrangian and the Lagrangian are related as,
\begin{align}
\hat{L}_\textrm{m} &= L +  h^4 \Phi_4(L) - h^4\Phi_2(\Phi_2(L))- \frac{h^4}{24} \ddot\theta_2(L)+ o(h^6), \label{equ:4mod}
\end{align}
where
\begin{align}
\Phi_2(L)  =  \theta_2(L)-\frac{1}{24}\ddot{L}, \quad \Phi_4(L)  = \theta_4(L) + \frac{7}{5760}{L}^{(4)}.  
\end{align}
Recall from equation (\ref{equ:modsurrsys}) that the evolution of the modified surrogate system differs from that of the original Euler-Lagrange equations in terms of order $h^4$,
Therefore, higher order derivatives of the system configuration ($\ddot{q}, q^{(3)}, \dots$) can be replaced in equation (\ref{equ:4mod}) resulting in a function of $q$ and $\dot{q}$ while maintaining explicit expressions for terms up to order $h^4$.
Note that $\Phi_1(L)$ and $\Phi_2(L)$ are linear operators such that $\Phi_i(\alpha L_1 + \beta L_2) = \alpha\Phi_i(L_1) + \beta\Phi_i(L_2)$. 
Furthermore, note that expressions (\ref{equ:mod2}) and (\ref{equ:4mod}) are equivalent.
However, fourth order terms are now explicitly given in equation (\ref{equ:4mod}). 
The surrogate Lagrangian, $\hat{L}$, is augmented such that these 4th order terms are accounted for and the \emph{fourth order surrogate} Lagrangian, $\hat{L}^{(4)}$, is derived, 
\begin{align}
\hat{L}^{(4)} &= \hat{L}  -  h^4 \Phi_4(L)+ h^4\Phi_2(\Phi_2(L))+ \frac{h^4}{24} \ddot\theta_2(L), \nonumber \\
 &= L -h^2 \Phi_2(L)  -  h^4 \Phi_4(L) + h^4\Phi_2(\Phi_2(L))+ \frac{h^4}{24} \ddot\theta_2(L).\label{equ:surr4lag}
\end{align}
Equation (\ref{equ:pre4mod}) is used to derive an expression for the \emph{ modified fourth order surrogate} Lagrangian, $\hat{L}^{(4)}_\textrm{m}$:
\begin{align} 
\hat{L}^{(4)}_\textrm{m} =& L + o(h^6).
\end{align}
The fourth order surrogate differential equation approximates the original differential equation up to order $h^6$ while the expression given in equation (\ref{equ:modsurrsys}) only does so up to order $h^4$. 
As before, the increase of accuracy was not achieved through an increase in complexity of the variational integrator. 
\emph{If convergence issues are ignored this process can be repeated indefinitely to produce integrators of arbitrary order.} 
However, at some point errors produced by, for example, floating point operations will become a limiting factor. 
The same process can also be applied to constrained and forced systems. 

\section{Numerical Experiments}

\subsection{Higher Order Surrogates}
The harmonic oscillator studied in Section \ref{sec:lowharmosc} is revisited to elucidate our discussion of higher order surrogate Lagrangians given in the previous section. 
As before, the resulting surrogate Lagrangian describes a mass-spring system with a different mass and spring constant. 
If the analysis shown in Section \ref{sec:higher} is repeated twice more an eighth order surrogate Lagrangian is obtained and parameterized as 
\begin{align}
M_{\textrm{s}} &= M  - K\frac{h^2}{12}- KM^{-1}K\frac{h^4}{720} -(KM^{-1})^2K\frac{h^6}{30240} - (KM^{-1})^3K\frac{h^8}{1209600},\\
K_{\textrm{s}} &= K +   KM^{-1}K\frac{h^2}{12} + (KM^{-1})^2K\frac{h^4}{120} + (KM^{-1})^3K\frac{17 h^6}{20160} + (KM^{-1})^4K\frac{31h^8}{362880}.
\end{align}
Note that surrogate Lagrangians of fourth and sixth order can be obtained by removing the appropriate terms in $M_{\textrm{s}}$  and $K_{\textrm{s}}$.
It is interesting to note that modifications of similar structure were report in \cite{mushtaq2014higher}.
However, the methodology presented there required the system to have a separable Lagrangian. 

\begin{figure}
    \centering
    \begin{subfigure}[]{0.5\columnwidth}
	\centering
        	\includegraphics[width=\columnwidth]{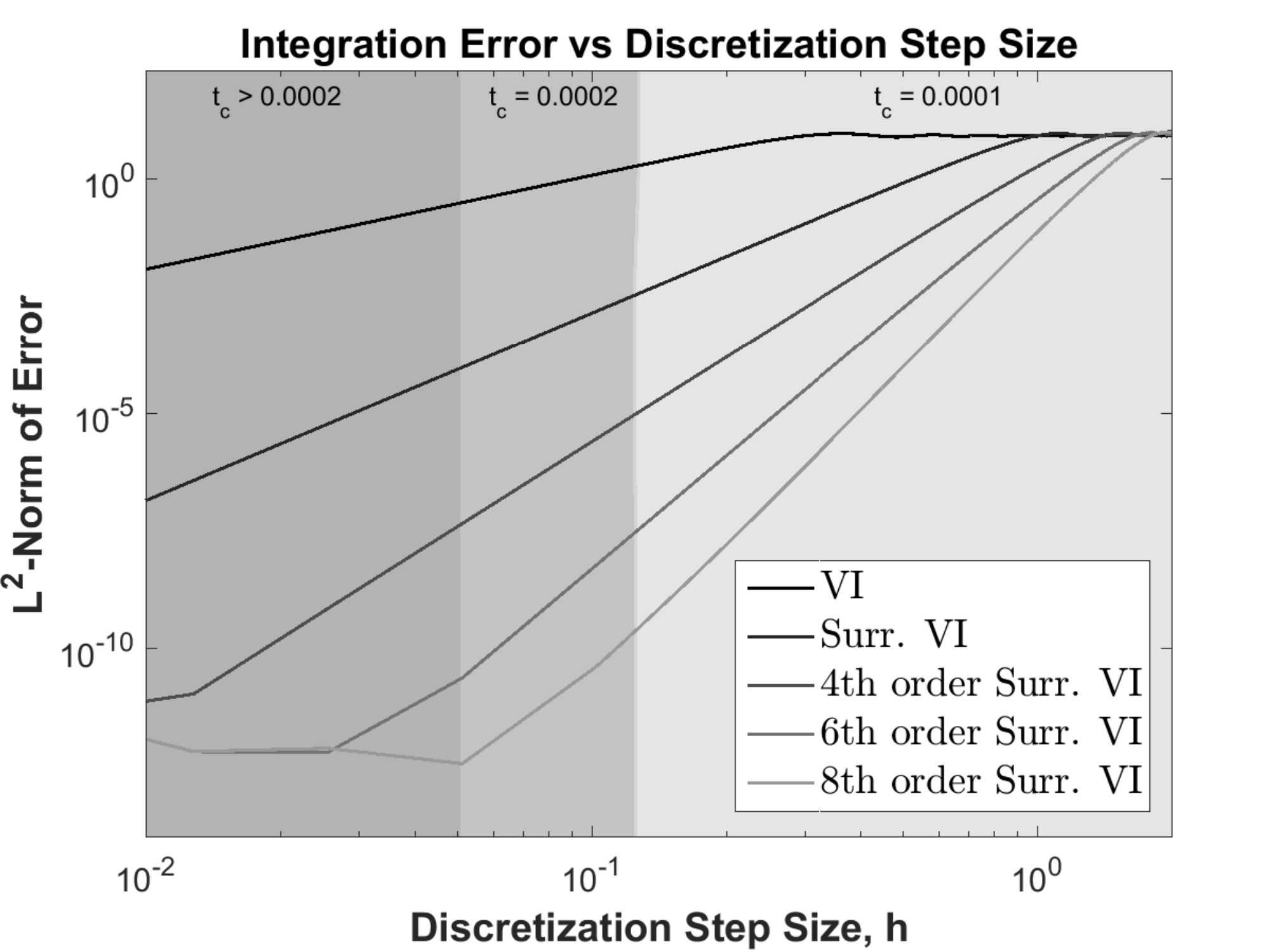}
	\caption{ }
    \end{subfigure}%
    \begin{subfigure}[]{0.5\columnwidth}
	\centering
        	\includegraphics[width=\columnwidth]{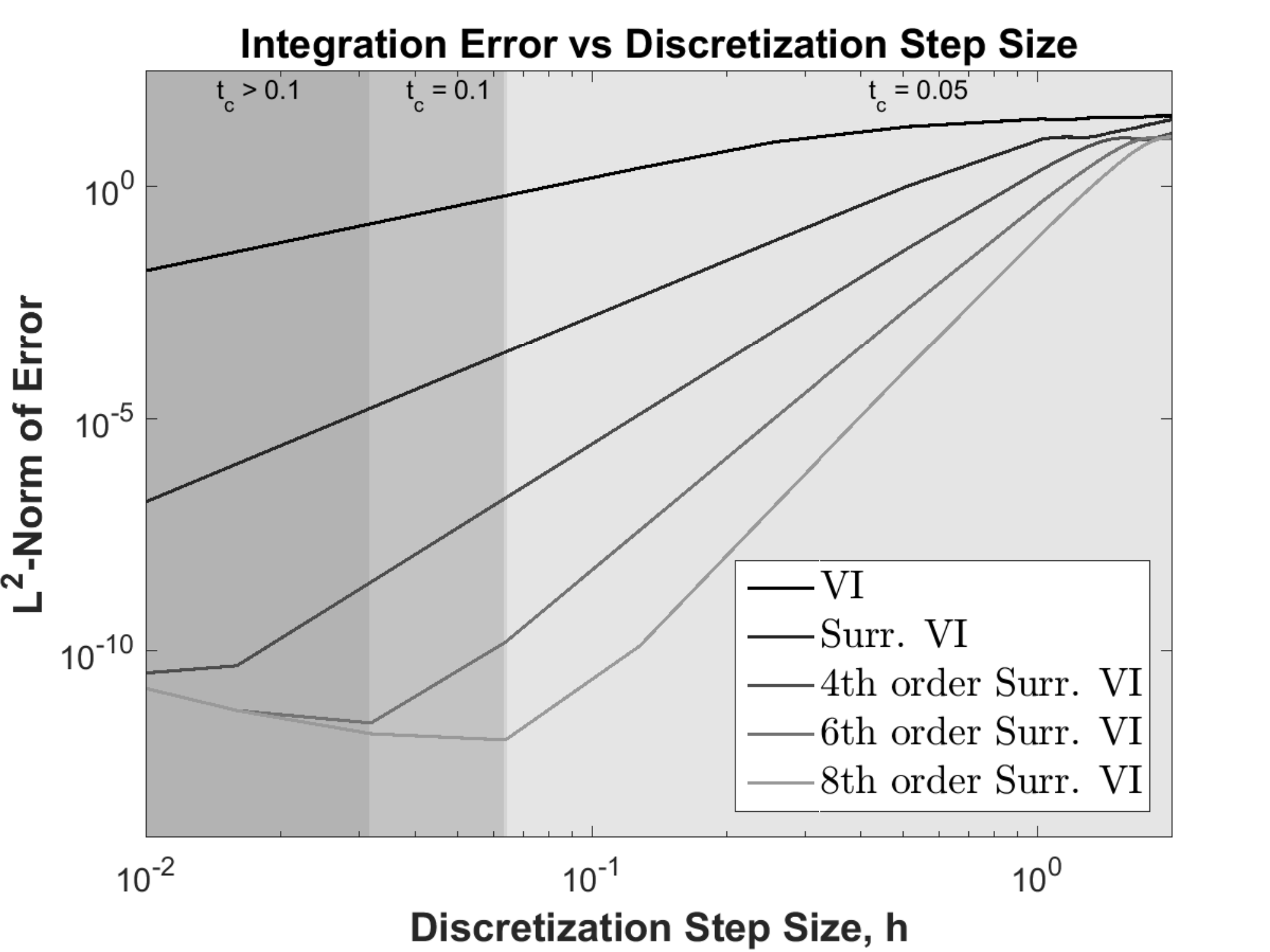}
	\caption{ }
    \end{subfigure}
    \caption{ The $L^2$-norm of the error as a function of the utilized discretization time step plotted on a contour of average computational times for (a) a one dimensional harmonic oscillator and (b) a four degrees of freedom mechanical system described by M and K in (\ref{4dofmk}), respectively. 
The nominal variational integrator and the surrogate variational integrators exhibit the expected order convergence of the $L^2$-norm of the error.
The increase in accuracy in the variational integrator is not associated with an added level of computational effort.
}
 \label{fig:higherspring}
\end{figure}

The initial condition of the system was set as $q(0)=0$ and $\dot{q}(0)=1$ and $M=1$ and $K=2$. 
The system was propagated for 150 seconds. 
Figure \ref{fig:higherspring} displays that convergence properties of the nominal variational integrator and 4 surrogate variational integrators. 
As before, the increase in accuracy is not accompanied by an added level of computational effort. 
Furthermore, each surrogate variational integrator displays the predicted order of convergence of the $L^2$-norm of the error. 

The benefits of the surrogate Lagrangian approach are not limited to single degree of freedom systems.
Consider a four degrees of freedom mechanical system described by mass and spring matrices given as
\begin{equation} \label{4dofmk}
M = 
        \begin{bmatrix}
                    2& 0.1 & 0 &  0.3 \\
 0.1&3&0.1&0\\
0&0.1&4.1&0.3\\  
0.3&0&0.3&4
        \end{bmatrix}, \quad K = \begin{bmatrix}
1&0.5&0&0.5\\ 
0.5&0.9&0.35&0\\ 
0&0.35&8.1&0.65\\  
0.5&0&0.65&2.1
        \end{bmatrix}\,
\end{equation}
such that its Lagrangian is given as ${L}=  \frac{1}{2}\dot q^\textrm{T}Mq - \frac{1}{2}q^\textrm{T}K q$. 
As before, the system was propagated for 150 seconds. 
Figure \ref{fig:higherspring} shows that the predicted orders of convergence where achieve by the implemented surrogate variational integrators. 

\subsection{Damped Harmonic Oscillator}
We now consider a damped mass-spring system with mass $M$, spring constant $K$, and damping coefficient $C$.
Its Lagrangian is given as 
\begin{eqnarray}
{L}&=  \frac{1}{2}M\dot q^2 - \frac{1}{2}K q^2,
\end{eqnarray}
and is subjected to a damping force given as
\begin{eqnarray}
F&=  -C\dot{q}.
\end{eqnarray}
Its surrogate Lagrangian is given as
\begin{eqnarray}
\hat{L}&=  \frac{1}{2}(M-\frac{h^2}{12}K+\frac{h^2}{12}\frac{C^2}{M})\dot q^2 - \frac{1}{2}(K+\frac{h^2}{12}\frac{K^2}{M})q^2,
\end{eqnarray}
and its surrogate discrete force, introduced in equation (\ref{equ:surroodiscreteforce}), is parameterized as 
\begin{eqnarray}
\hat{H}_1 = -b\dot{q} -\frac{h^2}{12}\frac{KC}{M}\dot{q}, \quad \hat{H}_2= \frac{h^2}{12}\frac{KC}{M}{q}.
\end{eqnarray}
As before, the surrogate Lagrangian describes  a mass-spring system with a  discretization time step dependent mass, spring constant, and damping coefficient.
The system was propagated using a variational integrator, a surrogate variational integrator, and a Runge-–Kutta fourth-order method and an analytic solution was used as a benchmark. 
The initial condition of the system was set as $q(0)=\sqrt{2}/2$ and $\dot{q}(0)=\sqrt{2}/2$ and $M=10$, $K=3$, and $C= 0.07$. 
The system was propagated for 300 seconds. 

\begin{figure}
    \centering
    \begin{subfigure}[]{0.5\columnwidth}
	\centering
        	\includegraphics[width=\columnwidth]{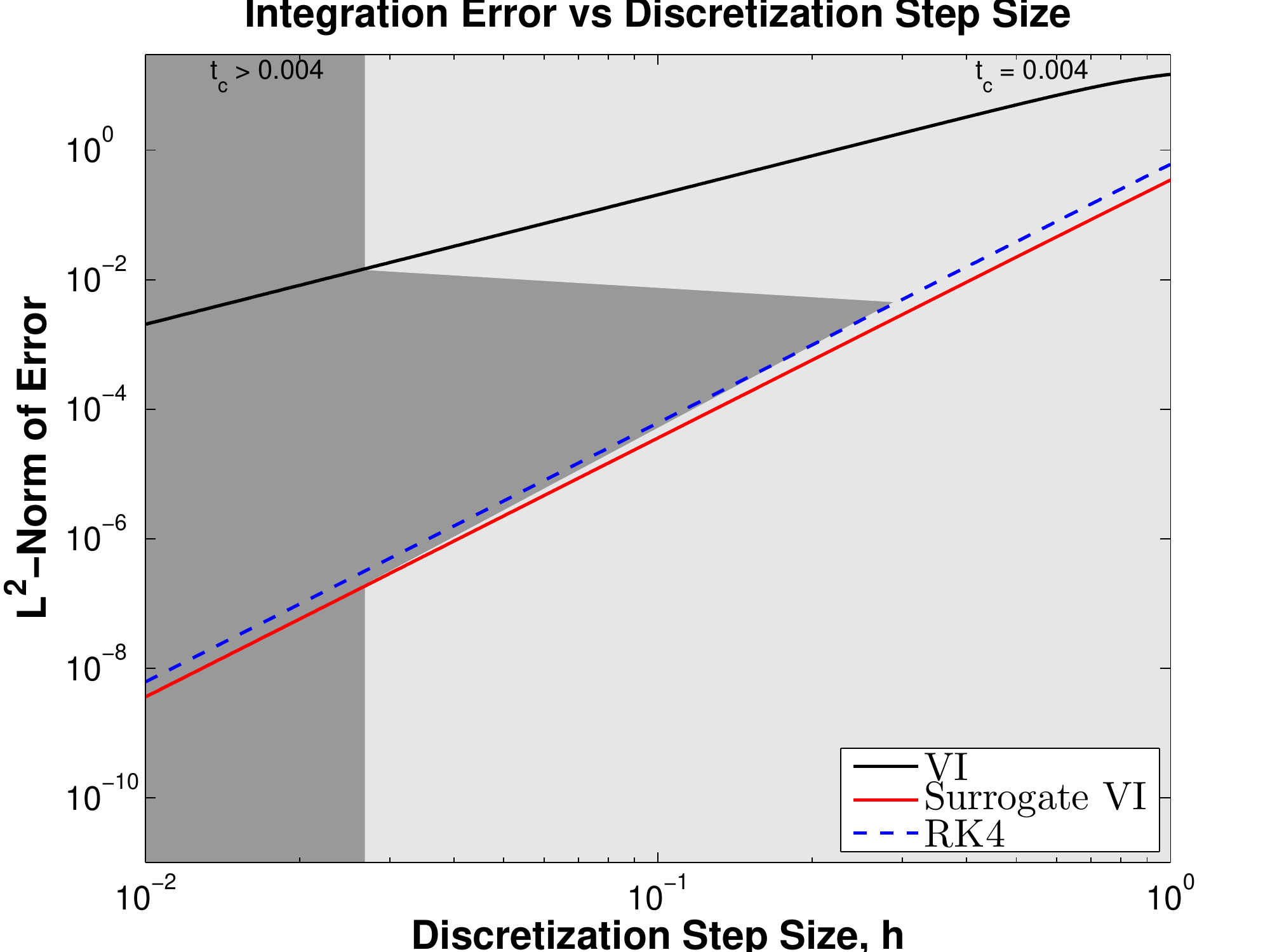}
	\caption{ }
    \end{subfigure}%
    \begin{subfigure}[]{0.5\columnwidth}
	\centering
        	\includegraphics[width=\columnwidth]{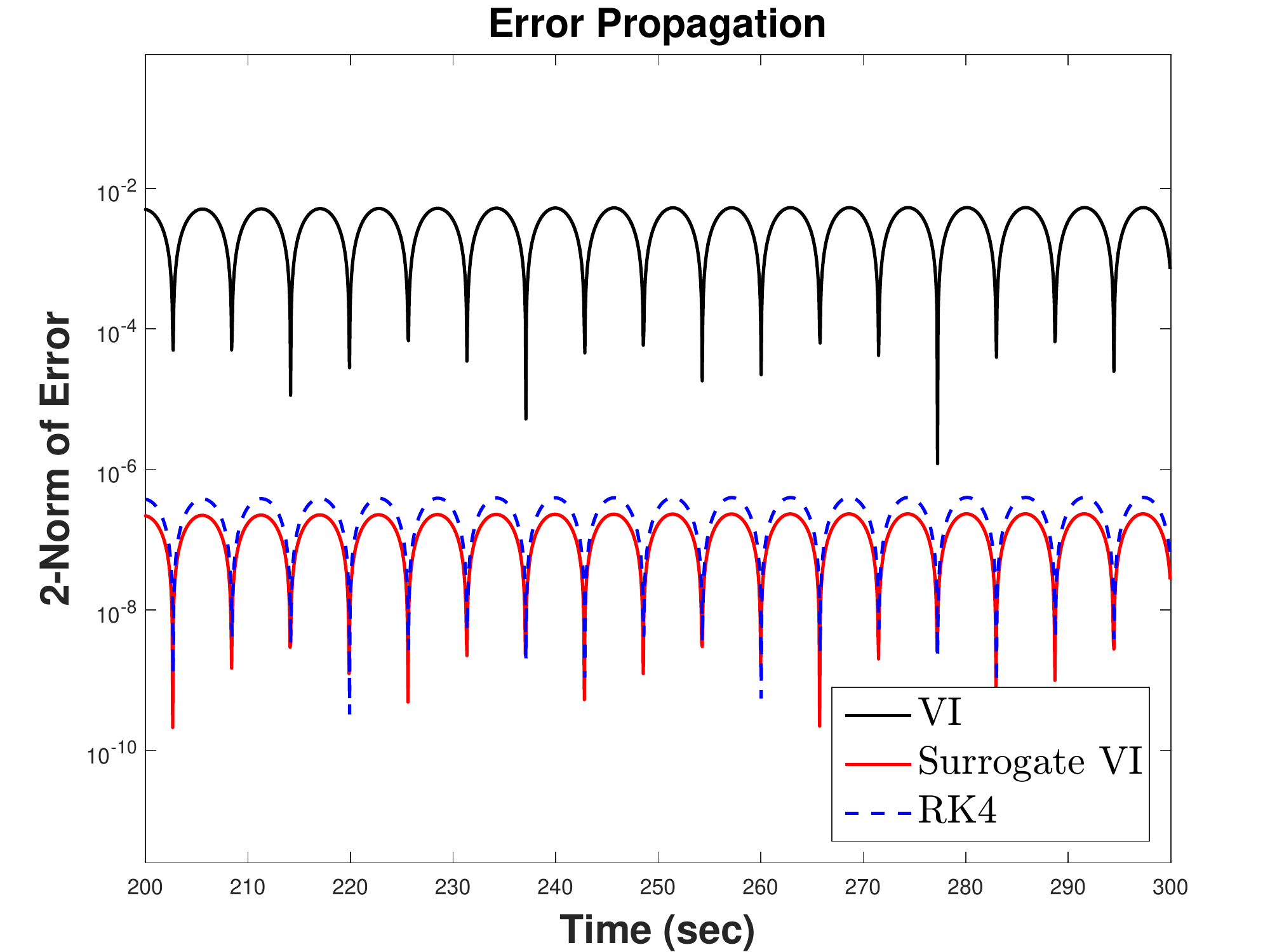}
	\caption{ }
    \end{subfigure}
    \caption{(a): The $L^2$-norm of the error as a function of the utilized discretization time step plotted on a contour of average computational times. 
The Runge--Kutta fourth-order method and the surrogate variational integrator exhibit a 4th order convergence of the $L^2$-norm of the error.
(b): The 2-norm of error for a particular execution of the numerical integrators when $h=0.05$ in the time interval $200\leq t \leq 300$.
}
    \label{fig:dampedspring}
\end{figure}

Figure \ref{fig:dampedspring} shows the predicted evolution of the system when $h=1.0$ and the $L^2$-norm of the error for a range of discretization time steps.
A contour of the average computational time for 100 executions of each integrator configuration is also displayed.
As before, the surrogate variational integrator exhibits a fourth order convergence and achieves the lowest $L^2$-norm in each contour depicted. 

Higher order convergence can be obtained by extending the procedure given in Section \ref{sec:higher} to forced systems. 
Therefore, a very accurate and computationally inexpensive simulation of a forced point mass is possible. 
Though point masses give an overly simplistic representation of most systems they nevertheless have many scientific and engineering applications. 
For example, finite element analysis methods will benefit from efficient and effective propagation of such systems \cite{FEM}.

\subsection{Pendulum: Cartesian Coordinates}
Single and double pendulums are considered to demonstrate the proposed methodology's utility for a nonlinear and, possibly, chaotic system. 
The position of the pendula are described with Cartesian coordinates and holonomic constraints are used to ensure pendulum lengths are maintained. 
The mass of each pendulum is assumed to be concentrated at the end of the link and is affected by a gravitational field.
The potential and kinetic energies of a single pendulum system are described as 
\begin{align}
V(q) = mgy, \quad T(\dot{q}) = \frac{1}{2}m(\dot{x}^2 + \dot{y}^2).
\end{align}
and its  holonomic constraint is given as
\begin{align}
c(q) = x^2 + y^2 - l.
\end{align}
Two half-explicit Runge-–Kutta integrators and the nominal and surrogate variational integrators were used to propagate the constrained systems\cite{arnold1998half}. 
In half-explicit Runge--Kutta methods the algebraic solution of the constraint variable is solve at each stage. 
Therefore, it is ensured that at each computed stage the given constraint is satisfied. 
However, the order of convergence of half-explicit methods do not, in general, correspond to the number of computed stages as would be expected in unconstrained Runge--Kutta methods.  
In fact, in order to obtain fourth-order convergence at least five stages are needed \cite{brasey1993half}.
In the proceeding examples a 4 stage half-explicit Runge--Kutta integrator that uses the standard fourth-order Butcher tableau and the half-explicit Runge--Kutta HEM4 algorithm, a fourth-order method with five stages developed in \cite{brasey1993half}, are evaluated alongside the variational integrators.

The benchmark to compute an integrator's $L^2$-norm of the error was the trajectory obtained by the respective integrator when $h= 1\times10^{-4}$.
The initial condition of the system was set as $q_1(t_0)=0, q_2(t_0)=\sqrt{l}, \dot q_1(t_0)=2$, and $\dot q_2(t_0)=0$ and system parameters were set as $m,l=1$ and $g=9.81$. 
The system was simulated for 10 seconds. 

\begin{figure}  
    \centering
    \begin{subfigure}[]{0.5\columnwidth}
	\centering
        	\includegraphics[width=\columnwidth]{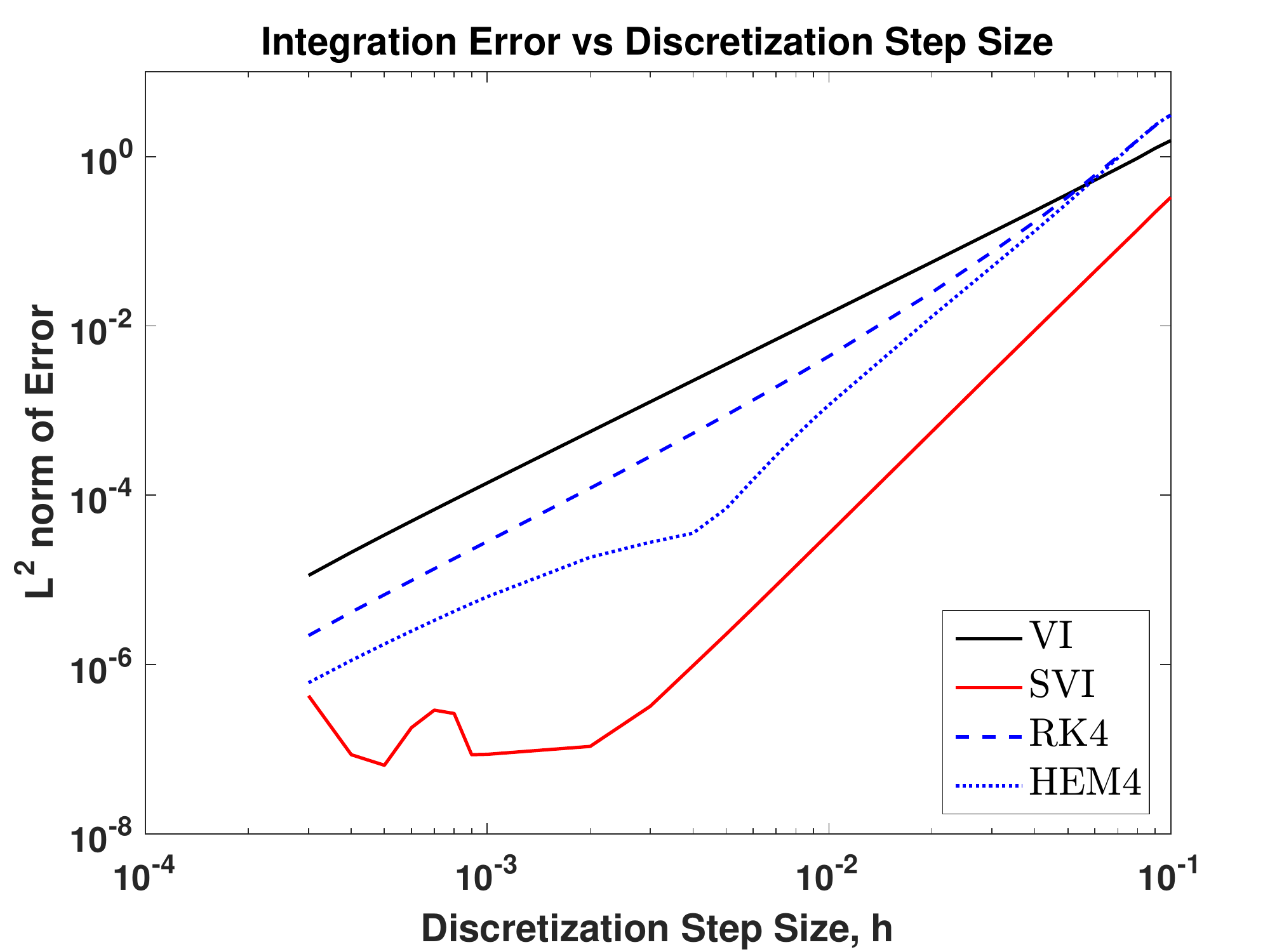}
	\caption{ }
    \end{subfigure}%
    \begin{subfigure}[]{0.5\columnwidth}
	\centering
        	\includegraphics[width=\columnwidth]{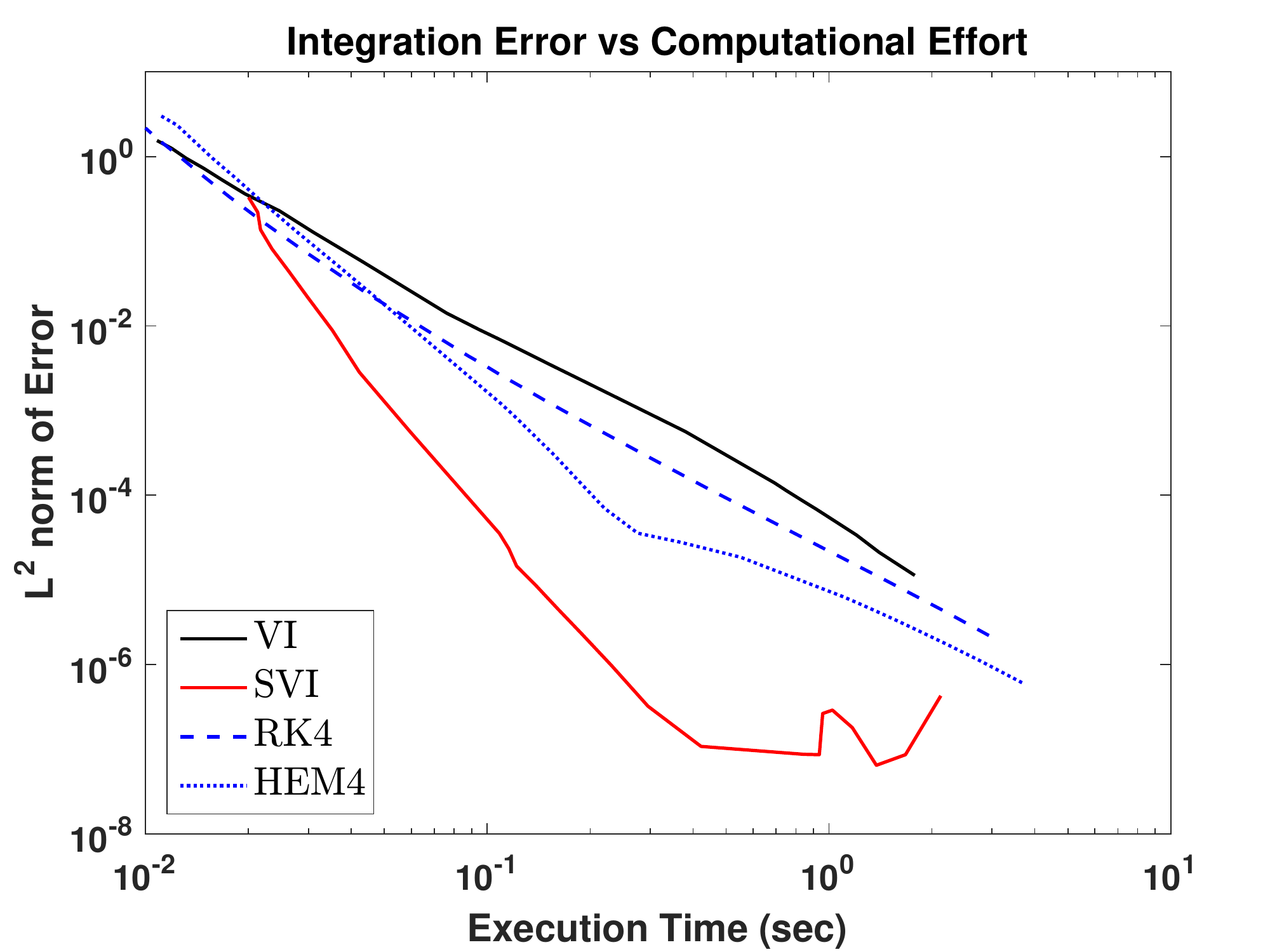}
	\caption{ }
    \end{subfigure}
    \caption{ (a): The $L^2$-norm of the error (benchmark trajectory obtained by the respective integrator when $h= 1\times10^{-4}$) as a function of the utilized discretization time step for a single pendulum. 
Note that the 4 stage half-explicit Runge--Kutta integrator displays second-order convergence, for the majority of the considered discretization step sizes, despite having fourth order complexity. 
The surrogate variational integrator obtains the smallest $L^2$-norm of the error at each discretization step size.
(b): The $L^2$-norm of the error as a function of the computational execution time for a single pendulum.
When the execution time exceeds $2.17\times10^{-2}$ seconds the surrogate variational integrator obtains the most accurate trajectory.
}
\label{fig:singpen}
\end{figure}

Figure \ref{fig:singpen} shows the $L^2$-norm of the error as a function of the discretization time step and the average computational time from 10 executions.
Note that the variational integrators retained the convergence properties seen earlier. 
However, as predicted in \cite{hairer1989numerical} the 4 stage half-explicit Runge--Kutta integrator displays second-order convergence, for the majority of the considered discretization step sizes, despite having fourth-order complexity. 
While the HEM4 algorithm achieves fourth-order convergence the surrogate variational integrator obtains a smaller $L^2$-norm of the error at each discretization step size.
Furthermore, when the execution time exceeds $2.17\times10^{-2}$ seconds the surrogate variational integrator obtains the most accurate trajectory.

\begin{figure}
    \centering
    \begin{subfigure}[]{0.5\columnwidth}
	\centering
        	\includegraphics[width=\columnwidth]{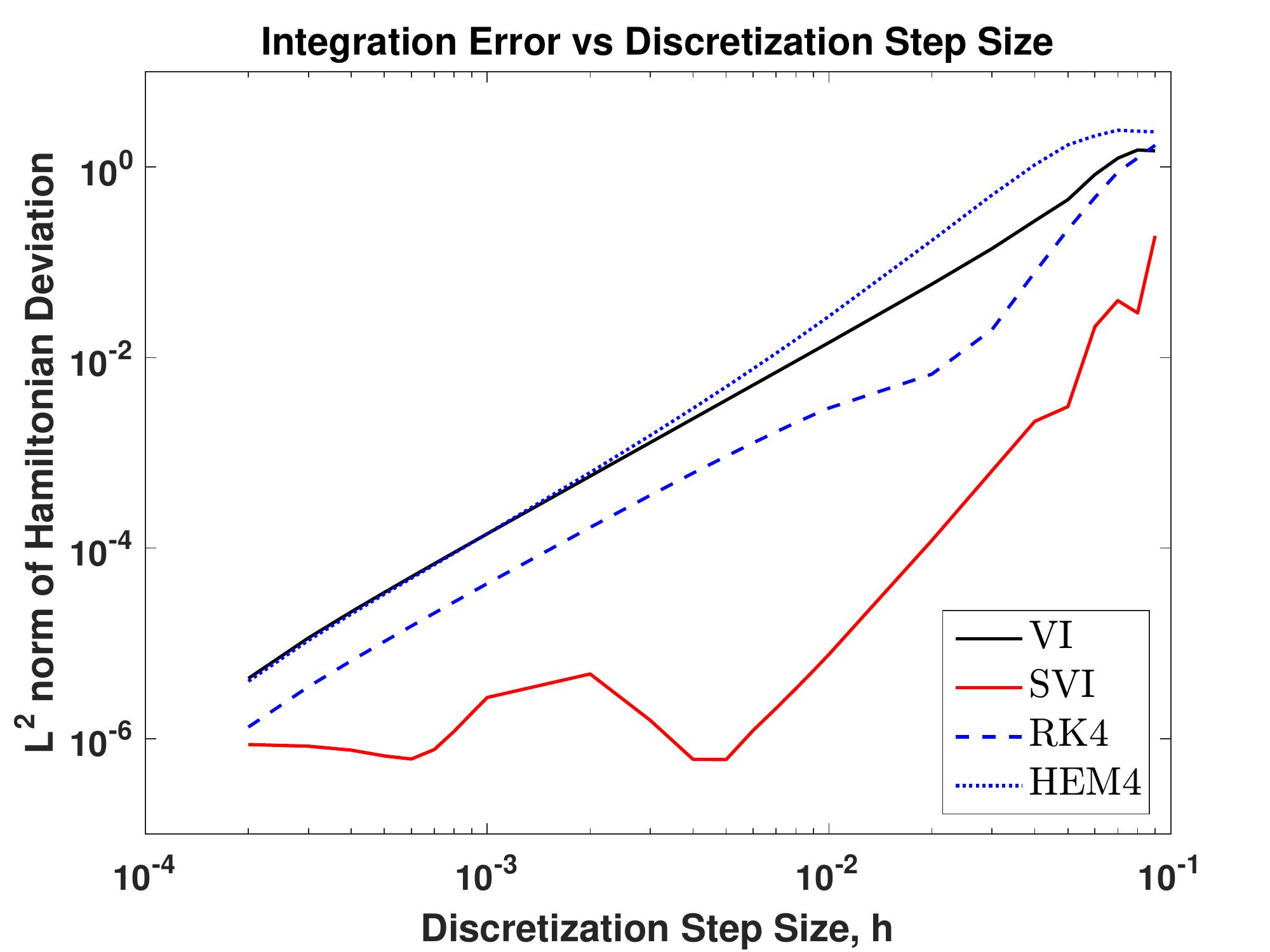}
	\caption{ }
    \end{subfigure}%
    \begin{subfigure}[]{0.5\columnwidth}
	\centering
        	\includegraphics[width=\columnwidth]{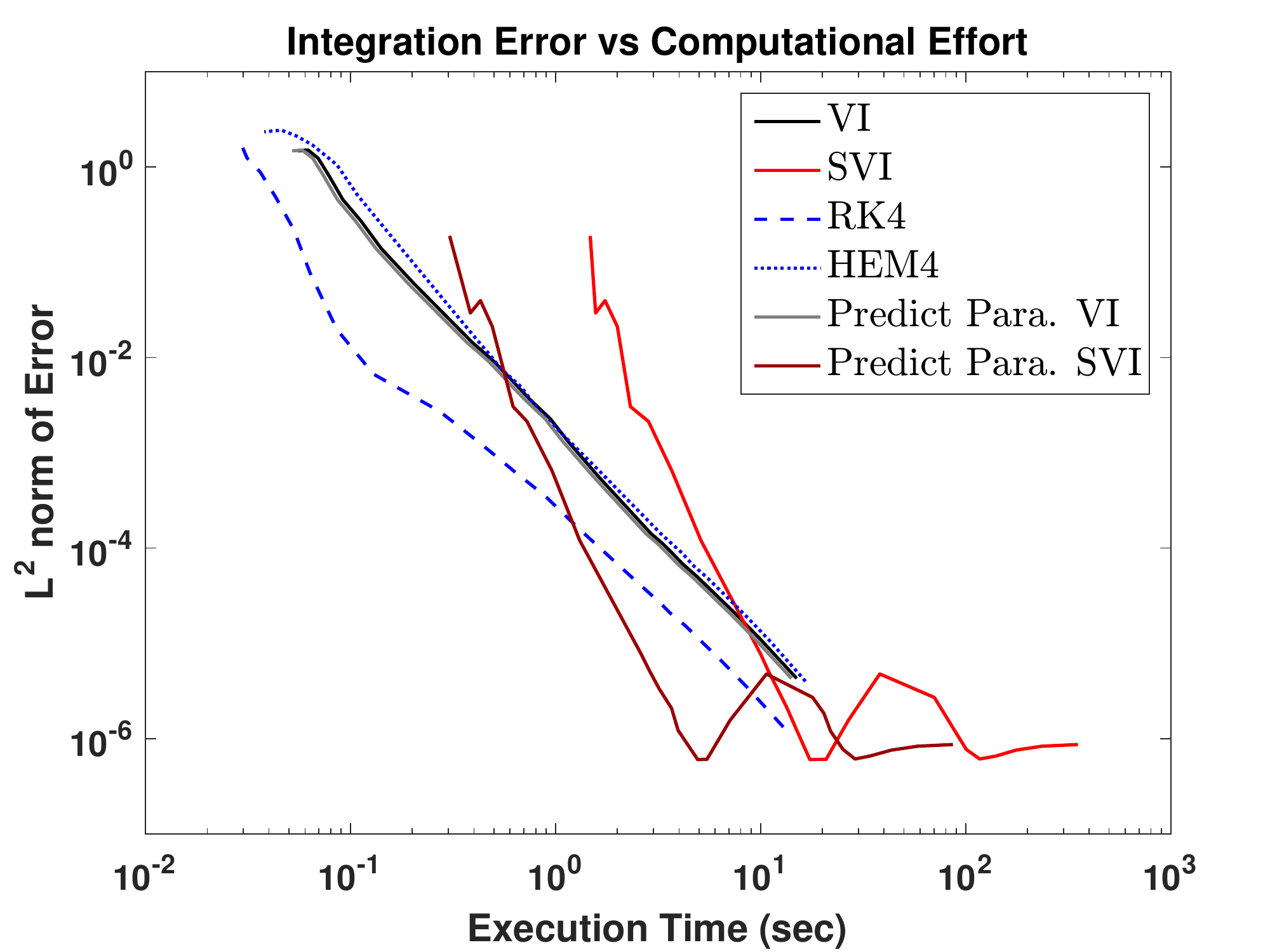}
	\caption{ }
    \end{subfigure}
    \caption{ (a): The $L^2$-norm of the error (benchmark trajectory obtained by the respective integrator when $h= 1\times10^{-4}$) as a function of the utilized discretization time step for a double pendulum. 
The variational integrators achieved their predicted rates of convergence.
However, the half-explicit methods display lower orders of convergence than predicted. 
(b): The $L^2$-norm of the error as a function of the computational execution time for a double pendulum.
The predicted performance of parallelized nominal and surrogate variational integrators are also shown.
It is supposed that parallelization reduces evaluation time  of $D_1L_\textrm{d}$ and $D_2D_1L_\textrm{d}$ by a factor of 5. 
}
\label{fig:doublepen2}
\end{figure}

In order to continue our examination the considered system was appended to include an additional pendulum (and holonomic constraint). 
The potential and kinetic energies of the appended system are described as 
\begin{align}
V(q) = m(gy_1 + gy_2), \quad T(\dot{q}) = \frac{1}{2}m(\dot{x}_1^2 + \dot{y}_1^2 + \dot{x}_2^2 + \dot{y}_2^2).
\end{align}
and its holonomic constraints are given as
\begin{align}
c(q) = [x_1^2 + y_1^2 - l, (x_1-x_2)^2 + (y_1-y_2)^2 - l]^\textrm{T}.
\end{align}
Again, the benchmark trajectories are obtained by the respective integrator when $h= 1\times10^{-4}$.
The initial condition of the system was set as $q_1(t_0),q_3(t_0)=0, q_2(t_0)=\sqrt{l}, q_2(t_0)=2\sqrt{l}, \dot q_1(t_0)=5$, and $\dot q_2(t_0),\dot q_3(t_0),\dot q_4(t_0)=0$ and system parameters were set as $m,l=1$ and $g=9.81$. 
The system was simulated for 10 seconds. 

Figure \ref{fig:doublepen2} shows the $L^2$-norm of the error as a function of the discretization time step and the average computational time from 10 executions.
While the variational integrators achieved their predicted rates of convergence the half-explicit methods are shown to under-perform. 
This suggest that the manner in which constraints are addressed by half-explicit methods may have a large effect on the accuracy of the propagated trajectory. 
Note that variational integrators simultaneously ensure that constraints and the preservation of mechanical energy are enforced when propagating the system's configuration. 
As before, the surrogate variational integrator provided the most accurate trajectory regardless of the discretization step size.  

Note that the surrogate variational integrator now requires more computational effort than the nominal variational integrator for any particular discretization time step. 
It is important to note that the increase in computational effort is not due to a change in the central integration scheme, but rather in an increase in the complexity of the evaluated terms. 
Specifically, the modification made to the Lagrangian cannot be described by a simple change in mass or a mechanical parameter. 
The resulting surrogate Lagrangian is a summation of differentiable functions (polynomials in this case) with more terms than the nominal Lagrangian.
Therefore, more computational effort is expended to evaluate the integration equation (\ref{back:fint}) and update the estimation of $q_{k+1}$.
However, the computational effort of the surrogate variational integrator can be mitigated if parallel computing is used to evaluate complex polynomials\cite{polyparra,gpu}.
Nevertheless, the surrogate variational integrator is the best choice, in terms of computational effort and accuracy, if an execution time of 12.9 seconds or more is acceptable. 

Table \ref{Table1} shows the amount of time each variational integrator requires to perform one step of integration ($(q_k,p_k)\rightarrow(q_{k+1},p_{k+1})$) and the time it spent evaluating $D_1L_\textrm{d}$ and $D_2D_1L_\textrm{d}$ when $h = 1\times10^{-3}$.
Note that since the variational integrators are implicit multiple evaluations of $D_1L_\textrm{d}$ and $D_2D_1L_\textrm{d}$ are needed. 
In this case, the Newton--Raphson method used $\epsilon_\textrm{tol} = 1\times10^{-9}$. 
Note that the evaluations of $D_1L_\textrm{d}$  and $D_2D_1L_\textrm{d}$ account for more than 94 percent of the surrogate variational integrator's computational effort. 
Therefore, parallelization of these operations will significantly reduce the computational effort of the surrogate variational integrator, but not that of the nominal variational integrator. 
Figure \ref{fig:doublepen2} and Table \ref{Table1} show the predicted computational effort if it is assumed that parallelization reduces $D_1L_\textrm{d}$ and $D_2D_1L_\textrm{d}$ evaluation time by a factor of 5. 
Implementation of a parallelized architecture in a graphics processing unit (GPU) can reduce the time needed to evaluate a sparse polynomial and its Jacobian by a factor of more than 10 \cite{gpu}. 
In this example the largest polynomial found in $\frac{\partial \hat{L}}{\partial q}$ and $\frac{\partial \hat{L}}{\partial \dot q}$ consist of approximately $290$ nominals. 
Therefore, the evaluated polynomials are quite sparse considering they are of degree 5 and contain 8 indeterminates (states of the system). 
Furthermore, instead of sequentially computing the 8 polynomials found in $\frac{\partial \hat{L}}{\partial q}$ and $\frac{\partial \hat{L}}{\partial \dot q}$ parallelization computation time can be further reduced by . 
Therefore, a hypothetical time reduction of a factor of 5 is reasonable. 
Note that in this case the surrogate variational is the best choice if an execution time of 1.21 seconds or more is acceptable. 
Therefore, parallelization of the surrogate variational integrator will greatly increase its utility. 

\begin{table}
\centering
   \begin{tabular}{ | c | c | c | c | c | c |}
   \hline
  Variational Integrator & Step & $D_1L_\textrm{d}$ & $D_2D_1L_\textrm{d}$ & (Predicted) Parallelized Step \\ \hline
     Nominal & $2.9\times10^{-4}$ & $1.3\times10^{-5}$ & $7.3\times10^{-6}$& $2.7\times10^{-4}$ \\ \hline
    Surrogate & $72\times10^{-4}$ & $16\times10^{-4}$ & $52\times10^{-4}$ & $18\times10^{-4}$ \\ 
    \hline
    \end{tabular}
\caption{Computational times of the nominal and surrogate variational integrators when $h = 1\times10^{-3}$. 
The Step column gives the amount of time each variational integrator requires to perform one step of integration ($(q_k,p_k)\rightarrow(q_{k+1},p_{k+1})$).
The $D_1L_\textrm{d}$ and $D_2D_1L_\textrm{d}$ columns give the total time the algorithm took evaluating the respective derivative (multiple evaluations needed).  
The Parallelized Step column give the predicted computational effort for a step if parallelization reduces $D_1L_\textrm{d}$ and $D_2D_1L_\textrm{d}$ evaluation time by a factor of 5. 
}
\label{Table1} 
\end{table}

\section{Conclusion}
We presented a methodology to increase the accuracy of variational integrators without adding complexity to their central integration scheme.
The methodology alters the variational integrator algorithm by replacing the system's Lagrangian with its surrogate Lagrangian.
Backward error analysis was used to formulate a general expression for surrogate Lagrangians.
Surrogate variational integrators exhibit fourth order of convergence despite remaining a second order integrator. 
The presented methodology can be used in a large class of dynamical systems, including those with external forcings and holonomic constraints, since no assumption on the Lagrangian's structure is made. 
Furthermore, an arbitrary order of convergence can be achieved by iteratively computing higher order surrogate Lagrangians. 

The presented numerical experiments demonstrate the efficacy of our approach. 
When compared to the nominal variational integrator, the surrogate variational integrator achieve higher order of convergence. 
When compared to the fourth-order Runge-Kutta method,  the surrogate variational integrator better preserved the system's mechanical energy.  
An increase in computational effort was seen in some experiements due to the complexity of evaluating the resulting surrogate Lagrangian and its derivatives. 
Future work will include parallelization of the methodology, analysis of computational complexity, investigation of connections to finite element analysis, and application to control and estimation algorithms. 

\section{Acknowledgements}
This work was supported by Army Research Office grant W911NF-14-1-0461.

\bibliographystyle{IEEEtran}
\bibliography{References}
\end{document}